\documentclass{article}

\author{Claire Levaillant\footnote{Part of this research was achieved during a stay of the author at
Institut Henri Poincar\'e in Paris. The author wishes to thank the
Institute for its great hospitality.}\\clairelevaillant@yahoo.fr}
\title{Reducibility of the Cohen--Wales
representation of the Artin group of type $D_n$}

\usepackage{amsmath, amssymb}
\usepackage{amsthm}
\usepackage{palatino}
\usepackage{amscd}
\usepackage{graphicx}
\usepackage{epsfig}
\usepackage{verbatim}

\newcommand{\mj}{\mathcal{J}}
\newcommand{\mx}{\mathcal{X}}
\newcommand{\mz}{\mathcal{Z}}
\newcommand{\my}{\mathcal{Y}}

\newcommand{\n}{\nu}
\newcommand{\Q}{\mathbb{Q}}

\newcommand{\be}{\beta}
\newcommand{\al}{\alpha}

\newcommand{\lra}{\longrightarrow}

\newcommand{\la}{\lambda}

\newcommand{\ovl}{\overline}
\newcommand{\unsurr}{\frac{1}{r}}
\newcommand{\unsur}{\frac{1}}

\newcommand{\nts}{\negthickspace}

\newcommand{\U}{\mathcal{U}}

\newcommand{\IH}{\mathcal{H}}

\newcommand{\lb}{\lbrace}
\newcommand{\rb}{\rbrace}
\newcommand{\noin}{\noindent}

\newcommand{\ih}{\IH_{F,r^2}}

\newcommand{\W}{\mathcal{W}}

\newcommand{\da}{\downarrow}

\newcommand{\cil}{\frac{n(n-3)}{2}}
\newcommand{\dbw}{\frac{(n-1)(n-2)}{2}}
\newcommand{\chl}{\frac{n(n-1)}{2}}

\newcommand{\cgwn}{\mathcal{C}\mathcal{G}\mathcal{W}(D_n)}
\newcommand{\wh}{\widehat}
\newcommand{\mcalh}{\mathcal{H}}

\newcommand{\mcalb}{\mathcal{B}}

\begin{document}
\maketitle \noin\textbf{Abstract.} Using knot theory, we construct a
linear representation of the CGW algebra of type $D_n$. This
representation has degree $n^2-n$, the number of positive roots of a
root system of type $D_n$. We show that the representation is
generically irreducible, but that when the parameters of the algebra
are related in a certain way, it becomes reducible. As a
representation of the Artin group of type $D_n$, this representation
is equivalent to the faithful linear representation of Cohen-Wales.
We give a reducibility criterion for this representation as well as
a conjecture on the semisimplicity of the CGW algebra of type $D_n$.
Our proof is computer-assisted using Mathematica.

\section{Introduction}
\subsection{Definitions and history}
In $2002$, Cohen and Wales showed the linearity of all the Artin
groups of finite type \cite{CW}. The same result was shown
independently by Digne in \cite{DI}. Linearity of a group means that
there exists a faithful linear representation of this group. In
other words, the group can be identified with a subgroup of
$GL_k(F)$ for some integer $k$ and some field $F$. If
$M=(m_{ij})_{1\leq i,j\leq m}$ is a Coxeter matrix of size $m$, the
Artin group of type $M$ is by definition the group with generators
$s_1,\,s_2,\,\dots,\,s_m$ and relations
$$\underset{\text{$m_{ij}$ terms}}{\underbrace{s_is_js_i\dots}}=\underset{\text{$m_{ij}$ terms}}{\underbrace{s_js_is_j\dots}}$$

\noindent The Artin group of type $A_{n-1}$ is the braid group $B_n$
on $n$ strands. In the past, several authors had tried to use the
$(n-1)$-dimensional Burau representation to show the linearity of
$B_n$. However, if this representation is faithful for $n=3$, it was
shown to be unfaithful for $n\geq 5$ (see \cite{MOO}, \cite{LP},
\cite{BI}). It is still unknown whether the Burau representation of
$B_4$ is faithful or not. Currently, the only known faithful linear
representation of the braid group $B_n$ for $n\geq 4$ is the
$\frac{n(n-1)}{2}$-dimensional Lawrence--Krammer representation.
This representation originated in the work of Ruth Lawrence in
\cite{RUT} and was later recovered by Daan Krammer in \cite{KR}.
Hence it is called the Lawrence--Krammer representation. The
Lawrence--Krammer space is a vector space spanned by vectors indexed
by the $\frac{n(n-1)}{2}$ positive roots of a root system of type
$A_{n-1}$. The beautiful result of linearity of the braid group
using the Lawrence--Krammer representation is due to Bigelow
\cite{BIG} and independently Krammer \cite{KR}. Soon after, Cohen
and Wales wanted to show that all the Artin groups of finite types
are linear groups. As part of their work in \cite{CW}, they
generalized the Lawrence--Krammer representation to types $D$ and
$E$. They then generalized Krammer's algebraic arguments to show the
faithfulness of these newly found representations. We will call
these representations the Cohen--Wales representations of respective
types $D$ and $E$. Cohen, Gijsbers and Wales shortly after in
\cite{CGW} build even more inequivalent representations of the Artin
groups of types $D$ and $E$. Except for the Cohen--Wales
representation, it is still unknown whether these representations
are faithful or not. These parameter-based representations that
include the Cohen--Wales representation all factor through the
Cohen-Gijsbers-Wales algebra (abbreviated CGW algebra), an algebra
that contains the Artin group. Working with generic parameters,
Cohen, Gijsbers and Wales show that these are all the irreducible
representations of a certain quotient of ideals of the CGW algebra.
In particular their work shows that the Cohen--Wales representation
is irreducible for generic parameters. The goal of the present paper
is to show that the Cohen--Wales representation of type $D_n$
becomes reducible when its two parameters are related in a certain
way. To do so, we use a knot theoretic approach to build a
representation of the CGW algebra of type $D_n$. As a representation
of the Artin group of type $D_n$, our representation is equivalent
to the Cohen--Wales representation. The fact that the two
representations are equivalent is a nontrivial fact that follows
from the results in \cite{CGW}. We use our representation to give a
reducibility criterion for the Cohen--Wales representation of type
$D_n$. Our work can be viewed as a generalization of \cite{CRAS},
where a reducibility criterion is given for the Lawrence--Krammer
representation of the braid group.

\subsection{Notations and main results}

The main result of this paper is the following.
\newtheorem*{maintheo}{Main Theorem}
\begin{maintheo}
Let $n$ be an integer with $n\geq 4$. Let $t$ and $r$ be two
non-zero complex numbers. \\Assume that
$\left|\begin{array}{l}r^{2k}\neq 1\;\text{for every integer}\;k
\;\text{with}\;1\leq k\leq n\\r^{2k}\neq -1\;\text{for every integer
$k$ with $1\leq k\leq n-1$}\end{array}\right.$ \\\\Then the faithful
Cohen--Wales representation of the Artin group of type $D_n$ based
on the parameters $t$ and $r$ is irreducible except when $$t\in\lb
r^{4n-4},r^{2n-4},-r^{2n-2},1,r^4,-1\rb,$$ when it is reducible.
\end{maintheo}

\noindent The restrictions on the parameter $r$ are natural ones and
we will explain them later on. They have a crucial role to play in
the paper.\\ 
We now introduce a few notations that relate to root systems of type
$D_n$. In what follows, $n$ is an integer with $n\geq 4$. The vector
space of the Cohen--Wales representation of type $D_n$ is spanned by
vectors indexed by the $n^2-n$ positive roots of a root system of
type $D_n$. We will number the nodes of the Dynkin diagram of type
$D_n$ as follows. We shall write $i\sim j$ if nodes $i$ and $j$ are
adjacent on the diagram. We denote by $r_1,\,\dots,\,r_n$ the simple
reflections.\vspace{-0.22cm}\begin{center}
\epsfig{file=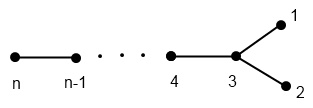,
height=2cm}\end{center}\vspace{-0.1cm} If $\al_1,\,\dots,\,\al_n$
denote the simple roots, then the positive roots are
\begin{itemize}
\item The $n$ simple roots $\al_1$, $\al_2$, $\al_3$, $\dots$, $\al_n$.
\item The ${n-1\choose 2}$ positive roots $\al_i+\dots\al_j$ with $2\leq
i<j\leq n$, of a root system of type $A_{n-1}$ on nodes
$2,\,\dots,\,n$.
\item The $(n-2)$ roots $\al_1+\al_3+\dots +\al_i$ with $i\geq 3$\\
The $(n-2)$ roots $\al_1+\al_2+\al_3+\dots+\al_i$ with $i\geq 3$.
\item The ${n-2\choose 2}$ roots
$\al_1+\al_2+2\,\al_3+\dots+2\,\al_i+\al_{i+1}+\dots+\al_j$ with
$3\leq i\leq n-1$ and $i+1\leq j\leq n$.
\end{itemize}
The spanning vectors of the Cohen-Wales space of type $D_n$ will be
denoted in the same order by $\widehat{w_{12}}$, $w_{12}$, $w_{23}$,
$\dots$, $w_{n-1,n}$ for the simple roots and by $w_{i-1,j}$,
$\widehat{w_{1,i}}$, $\widehat{w_{2,i}}$, $\widehat{w_{i,j}}$ for
the other positive roots. The significance of the indices will be
explained in detail later. For now, notice that a spanning vector
carries a hat if and only if $1$ is in the support of the positive
root, that is the coefficient of $\al_1$ in the positive root is
nonzero. \\
Let $g_1,\,g_2,\,g_3,\,\dots,\,g_n$ be the generators of the Artin
group $\mathcal{A}(D_n)$ of type $D_n$. As read on the Dynkin
diagram, the defining relations are as follows.
$$(A)\left\lbrace\begin{array}{l}\text{The $g_i$'s with $2\leq i\leq n$ satisfy the
braid relations}\\ \begin{array}{cccc}g_1\,g_3\,g_1&=&g_3\,g_1\,g_3&\\
g_1\,g_k&=&g_k\,g_1&\text{when $k\neq
3$}\end{array}\end{array}\right.$$ \noindent The CGW algebra $\cgwn$
of type $D_n$ is an algebra with two parameters $l$ and $m$ that
contains $\mathcal{A}(D_n)$, and contains other elements
$e_1,\,e_2,\,e_3,\,\dots,\,e_n$. These elements are related to the
generators $g_i$'s by $$m\,e_i=l\,(g_i^2+m\,g_i-1)$$ The other
defining relations of the algebra are as follows.
$$(DL)\left\lbrace\begin{array}{cccc}
g_i\,e_i&=&\unsur{l}\,e_i&\text{for all $i$}\\
e_ig_je_i&=&l\,e_i&\;\text{when $i\sim j$}
\end{array}\right.$$
Some selected immediate consequences of these definitions are the
following (see \cite{CGW}).
$$\left\lbrace\begin{array}{cccc}
e_i\,g_i&=&\nts\nts\nts\nts\unsur{l}\,e_i&\text{for all $i$}\\
g_i-g_i^{-1}&=&\nts\nts\nts m\,(e_i-1)&\text{for all $i$}\\
e_ie_je_i&=&\nts\nts\nts e_i&\text{when $i\sim j$}\\
g_jg_ie_j&=&e_ig_jg_i\;\;=\;\;e_ie_j&\text{when $i\sim j$}\\
e_je_ig_j&=&e_jg_i^{-1}\;\;=\;\;e_jg_i+m\,(e_j-e_je_i)&\text{when $i\sim j$}\\
g_je_ie_j&=&g_i^{-1}e_j\;\;=\;\;g_ie_j+m\,(e_j-e_ie_j)&\text{when $i\sim j$}\\
e_i^2&=&\nts\nts\nts\delta\,e_i&\text{with
$\delta=1-\frac{l-\unsur{l}}{m}$}
\end{array}\right.$$

Informations that relate to rank or cellularity can be found in
\cite{BAL2}. The CGW algebra of type $D_n$ is a generalization of
the BMW algebra to type $D_n$. The BMW algebra is named after
Birman, Murakami and Wenzl. It was introduced by Birman and Wenzl in
\cite{BW} and independently by Murakami in \cite{MUR}. It has the
same defining relations as above except they must be read on a
Dynkin diagram of type $A_{n-1}$. The BMW algebra is in connection
with a polynomial link invariant, namely the Kauffman polynomial. An
important feature of the Kauffman polynomial is that it can
distinguish oriented links that the other polynomials can't
distinguish. Birman and Wenzl wanted to build an algebra equipped
with a trace so that the Kauffman polynomial \cite{KAUF} is after
appropriate renormalization that trace, in the same way the
$2$-variable generalization of the Jones polynomial \cite{JONES},
namely the HOMFLY polynomial \cite{HOMFLY} was after renormalization
the trace on the Hecke algebra of the symmetric group. If we define
the Hecke algebra $\mathcal{H}(D_n)$ as the algebra with generators
$g_1,\,\dots,\,g_n$ that satisfy the Artin relations $(A)$ above and
the relations $g_i^2+m\,g_1=1$ for all $i$, we note that a quotient
of the CGW algebra is the Hecke algebra. This quotient of $\cgwn$
will play a critical role when studying the reducibility of the
representation of $\cgwn$ that we build. \noindent We next give the
expression of this representation. Its construction will be
explained in detail in $\S 2.2$. We introduce a new indeterminate
$r$ that is related to $m$ by the relation $m=\unsurr-r$. When later
specialized to non-zero complex numbers, $r$ and $-\unsurr$ will
hence be the two complex roots of the polynomial $X^2+m\,X-1$. We
choose $\Q(l,r)$ as our base field. As $r^2+m\,r=1$, the field
$\Q(l,r)$ is then a left $\mcalh(D_n)$-module for the action given
by $g_i.1=r$. We will denote by $V_n$ the Cohen-Wales space of type
$D_n$. The space $V_n$ is the vector space spanned over $\Q(l,r)$ by
the vectors $w_{ij}$'s and $\widehat{w_{ij}}$'s with $1\leq i<j\leq
n$. It will be convenient to introduce the following notation: by
$\ovl{w_{st}}$ for some integers $s$ and $t$ with $1\leq s<t\leq n$,
we mean that $\ovl{w_{st}}$ is either $w_{st}$ or $\wh{w_{st}}$. For
more clarity in the writing, we also sometimes add a coma between
the two indices $s$ and $t$.

\newtheorem{theo}{Theorem}
\begin{theo}
The following map $\n^{(n)}$ $$\begin{array}{ccc}
\begin{array}{l}
\cgwn\\
\\
g_i
\end{array}&\begin{array}{l} \stackrel{\n^{(n)}}{\lra}\\\\\longmapsto\end{array}&\begin{array}{l}
\text{End}_{\Q(l,r)}(V_n)\\\\\n_i
\end{array}\end{array}$$
defines a representation of the CGW algebra of type $D_n$ in the
Cohen-Wales space of type $D_n$. The actions are given as follows.
First, the action by $g_1$ is special and needs to be formulated
apart.
\begin{eqnarray}
\forall 3\leq i<j\leq n,\;\qquad\n_1(\wh{w_{ij}})&=&\begin{split}
m\,r^{j-5}(\wh{w_{1,i}}-w_{1,i})
\\+m\,r^{j-4}(\wh{w_{2,i}}-w_{2,i})\\+m\,r^{i-3}(\wh{w_{1,j}}-w_{1,j})\\+m\,r^{i-2}(\wh{w_{2,j}}-w_{2,j})\\
+m^2\,(r^{i+j-8}+r^{i+j-6})(\wh{w_{12}}-w_{12})\\+r\,\wh{w_{ij}}
\end{split}\\
\forall j\geq 3,\;\qquad\n_1(\wh{w_{1j}})&=&w_{2j}\\
\forall j\geq
3,\;\qquad\n_1(\wh{w_{2j}})&=&w_{1j}+m\,r^{j-3}\,w_{12}-m\,w_{2j}\\
\forall j\geq
3,\;\qquad\n_1(w_{2j})&=&\wh{w_{1j}}+m\,r^{j-3}\,\wh{w_{12}}-m\,w_{2j}\\
\forall j\geq
3,\;\qquad\n_1(w_{1j})&=&\begin{split}\wh{w_{2j}}+m\,r^{j-4}(\wh{w_{12}}-w_{12})\\+m\,(\wh{w_{1j}}-w_{1j})\end{split}\\
\n_1(\wh{w_{12}})&=&\unsur{l}\,\wh{w_{12}}\\
\text{In all the other
cases,}\;\n_1(\overline{w_{s,t}})&=&r\,\overline{w_{s,t}}
\end{eqnarray}
Second, the action by $g_2,\,g_3,\,\dots,\,g_n$ is determined by the
following expressions.
\begin{eqnarray}
\forall\,t\geq i+2,\;\qquad\n_{i+1}(\ovl{w_{i,t}})&=&\ovl{w_{i+1,t}}\\
\n_{j+1}(\ovl{w_{s,j}})&=&\ovl{w_{s,j+1}}\\
\n_{i+1}(\wh{w_{i,i+1}})&=&r\,\wh{w_{i,i+1}}\\
\n_{i+1}(w_{i,i+1})&=&\unsur{l}\,w_{i,i+1}\\
\n_i(\wh{w_{i,t}})&=&\wh{w_{i-1,t}}+\frac{m}{l}\,r^{t-i-2}\,w_{i-1,i}-m\,\wh{w_{i,t}}\\
\forall\,s\leq j-2,\;\qquad\n_j(\wh{w_{s,j}})&=&\wh{w_{s,j-1}}+\frac{m}{l\,r^{j-s-2}}\,w_{j-1,j}-m\,\wh{w_{s,j}}\\
\n_i(w_{i,t})&=&w_{i-1,t}+m\,r^{t-i-1}\,w_{i-1,i}-m\,w_{i,t}\\
\forall\,s\leq
j-2,\;\qquad\n_j(w_{s,j})&=&w_{s,j-1}+\frac{m}{l\,r^{j-s-2}}\,w_{j-1,j}-m\,w_{s,j}\\
\text{In all the other
cases,}\;\n_k(\ovl{w_{s,t}})&=&r\,\ovl{w_{s,t}}
\end{eqnarray}
As a representation of the Artin group, up to the change of
parameters $l=r^3\,t^{-1}$ and up to some rescaling of the
generators, this representation is equivalent to the Cohen-Wales
representation with parameters $t$ and $r$ that was built and used
in \cite{CW} to show the linearity of the Artin group of type $D_n$.
\end{theo}
\noindent In what follows, $\ih(n)$ denotes the Iwahori Hecke
algebra of the symmetric group $Sym(n)$ with parameter $r^2$ over
the field $\Q(l,r)$ as defined in \cite{MAT}. The following theorem
gives a reducibility criterion for the above representation under
some assumption of semisimplicity for the Hecke algebras $\ih(n)$
and $\mcalh(D_n)$.
\begin{theo}
Let $n$ be an integer with $n\geq 4$. Let $l$, $m$ and $r$ be three
non-zero complex numbers with $m=\unsurr-r$. \\
$(i)$ Assume that the Hecke algebras $\mcalh(D_n)$ and $\ih(n)$ are
semisimple. So assume that $r^{2k}\neq 1$ for every integer $k$ with
$1\leq k\leq n$ and $r^{2k}\neq -1$ for every integer $k$ with
$1\leq k\leq n-1$. Then, $\n^{(n)}$ is irreducible except when
$$l\in\bigg\lb\unsur{r^{4n-7}},\unsur{r^{2n-7}},-\unsur{r^{2n-5}},r^3,\unsurr,-r^3\bigg\rb,$$
when it is reducible. \\$(ii)$ For these values of the parameters
and the values for which $r$ has been replaced by $-\unsurr$, the
CGW algebra $\cgwn$ of type $D_n$ of \cite{CGW} with parameters $l$
and $m$ over the field $\Q(l,m)$ is not semisimple.
\end{theo}
\noindent Note Theorem $1$ together with point $(i)$ of Theorem $2$
imply the Main Theorem. On the way, we further show the following
theorems on the dimensions.\\\\ \textbf{Key assumption: until the
end of the paper, we assume that the Hecke algebras
$\mathbf{\mcalh(D_n)}$ and $\mathbf{\ih(n)}$ are semisimple}.

\begin{theo}\textbf{(Existence of a one-dimensional invariant
subspace)} Let $n$ be an integer with $n\geq 4$. In $V_n$ there
exists a one-dimensional invariant subspace if and only if
$l=\unsur{r^{4n-7}}$. If so, it is unique and it is spanned over
$\Q(l,r)$ by the vector
$$u=\sum_{1\leq i<j\leq n}r^{i+j}\,\big(\wh{w_{ij}}+r^{2n-4}\,w_{ij}\big)$$
\end{theo}

\begin{theo}\textbf{(Existence of an irreducible $(n-1)$-dimensional
invariant subspace)} Let $n$ be an integer with $n\geq 5$. In $V_n$
there exists an irreducible $(n-1)$-dimensional invariant subspace
if and only if $l=\unsur{r^{2n-7}}$. If so it is unique and it is
spanned over $\Q(l,r)$ by the vectors $v_i$'s, $1\leq i\leq n-1$
with
\begin{align*}v_i=(r^{2n-6}-\unsur{r^2})\,w_{i,i+1}&+\sum_{j=i+2}^n
r^{j-i-4}\bigg\lb
(w_{i+1,j}-r\,w_{i,j})+r^2\,(\wh{w_{i+1,j}}-r\,\wh{w_{i,j}})\bigg\rb\\
&+\sum_{s=1}^{i-1}r^{s-i}\bigg\lb
r^{2n-6}(w_{s,i+1}-r\,w_{s,i})+(\wh{w_{s,i+1}}-r\,\wh{w_{s,i}})\bigg\rb\end{align*}
\end{theo}

\begin{theo}\textbf{(Existence of an irreducible $n$-dimensional
invariant subspace)}\\$(i)$ Let $n$ be an integer with $n\geq 4$ and
$n\neq 5$. If there exists an irreducible $n$-dimensional
invariant subspace inside $V_n$, then $l=-\unsur{r^{2n-5}}$.\\
$(ii)$ \textit{(Case $n=5$)} If there exists an irreducible
$5$-dimensional invariant subspace inside $V_5$, then $l\in\lb
-\unsur{r^5},r^3\rb$. If so, it is unique.
\end{theo}

\begin{theo}\textbf{(Existence of an irreducible $\cil$-dimensional
invariant subspace)}. Let $n$ be an integer with $n\geq 6$. If there
exists an irreducible $\cil$-dimensional invariant subspace inside
$V_n$ then $l=r^3$.
\end{theo}

\begin{theo}\textbf{(Existence of an irreducible
$\frac{n(n-1)}{2}$-dimensional invariant subspace)}. Let $n$ be an
integer with $n\geq 5$. If $l=\unsurr$, there exists a unique
irreducible $\frac{n(n-1)}{2}$-dimensional invariant subspace inside
$V_n$. Moreover, it is spanned over $\Q(l,r)$ by the vectors
$t_{ij}=w_{ij}-\wh{w_{ij}}$, $1\leq i<j\leq n$.
\end{theo}

\noindent The theorems above have been stated in increasing order
for the dimensions when $n\geq 5$. While Theorems $3$ and $4$
provide necessary and sufficient conditions for the existence,
Theorems $5,6$ and Theorem $7$ only provide a necessary condition
and a sufficient condition respectively. We show along the proof of
Theorem $2$ point $(i)$ that when the representation is reducible,
the action on a proper invariant subspace is a $\mcalh(D_n)$-action.
When $\mcalh(D_n)$ is semisimple, which we assume in this paper, the
irreducible representations of $\mcalh(D_n)$ are indexed by
unordered double partitions $(\la,\mu)$ of $n$ as in Theorem $1.5$
of \cite{Hu1}. Their degrees is given by the number of standard
Young tableaux of shape $(\la,\mu)$. By standard, we mean that the
integers ranging from $1$ to $n$ must be filled in the tableau with
the numbers increasing along the rows and down the columns. $(0)$
denotes the empty partition. The classes of irreducible
$\mcalh(D_n)$-modules are called Specht modules. In \cite{Clas}, we
give the complete classification of the invariant subspaces of the
representation in terms of Specht modules.\\

We end this introduction by presenting a conjecture that relates to
point $(ii)$ of Theorem $2$. It gives a semisimplicity criterion for
the CGW algebra of type $D_n$ in the same spirit as existing
criteria for type $A$. Let's briefly recall these criteria in type
$A$. In \cite{WEN}, Hans Wenzl was the first to discuss the
semisimplicity of the Birman--Murakami--Wenzl algebra. He considers
the BMW algebra with nonzero complex parameters $l$ and $m$. He
shows the following result.
\newtheorem*{thm}{Theorem}
\begin{thm}\textbf{[Wenzl], $1988$}
The BMW algebra with nonzero complex parameters $l$ and $m$ is
always semisimple except possibly if $r$ is a root of unity of if
$l$ is some power of $r$, where $r$ is a complex root of the
polynomial $X^2+m\,X-1$. \end{thm} \noin Some of these powers are
identified twenty years later in the Ph.D. thesis of \cite{THE} and
Theorem $2$, point $(ii)$ of the present paper can be viewed as a
generalization of Theorem $2$ of \cite{CLDBW}. We recall below this
result in type $A$.
\begin{thm}\textbf{[Levaillant-Wales], $2008$}
Let $n$ be an integer with $n\geq 3$. Let $m$, $l$ and $r$ be three
nonzero complex numbers with $m=\unsurr-r$.\\
\\$1)$ Suppose $n\geq 4$. If $r^{2k}=1$ for some $k\in\lb
2,\dots,n\rb$ or if $l$ belongs to the set of values
$r,-r^3,\unsur{r^{2n-3}},\unsur{r^{n-3}},-\unsur{r^{n-3}},-r^{2n-3},r^{n-3},-r^{n-3},\unsur{r^3},-\unsurr$,
the BMW algebra $BMW_n$ of type $A_{n-1}$ with parameters $l$ and
$m$ over the field $\Q(l,r)$ is not semisimple.\\\\
$2)$ If $r^4=1$ or $r^6=1$ or if $l\in\lb -r^3,\unsur{r^3},1,-1\rb$,
the algebra $BMW_3$ with parameters $l$ and $m$ over the field
$\Q(l,r)$ is not semisimple.
\end{thm}
\noin Simultaneously and independently, using cellularity techniques
that were first introduced by John Graham and Gus Lehrer in
\cite{GRLR}, Hebing Rui and Mei Si find a complete criterion of
semisimplicity for the BMW algebra. Their work is based on the
groundbreaking work of John Enyang \cite{ENY}, where he constructs a
cellular basis for the BMW algebra. Let's recall here Theorem $B$ of
\cite{RUI}
\begin{thm}\textbf{[Rui-Si], $2009$} Let $n$ be an integer with $n\geq 3$. Let $\mcalb_n$ be the
Birman-Murakami-Wenzl with parameters $l$ and $r$. \\\\
$a)$ Suppose $l\not\in\lb -r,\unsurr\rb$. Then $\mcalb_n$ is
semisimple if and only if $o(r)>2n$ and
$$l\not\in\bigcup_{k=3}^{n}\lb r^{3-2k},\pm
r^{3-k},-r^{2k-3},\pm r^{k-3}\rb$$.\\
$b)$ Assume $l\in\lb -r,\unsurr\rb$. Then,\\\\ $\mcalb_n$ is not
semisimple if $n$ is either even or odd with $n\geq 7$.\\
$\mcalb_3$ is semisimple if and only if $o(r)>6$ and $r^4\neq -1$\\
$\mcalb_5$ is semisimple if and only if $o(r)>10$ and $r^6\neq -1$
and $r^8\neq -1$.
\end{thm}

\noin In their Theorem, case $b)$ is different and is the case when
$e_i^2=0$. \\
Based on Theorem $2$ of this paper, and in the spirit of the
Theorems [Wenzl] and [Rui-Si] stated above, we give the following
conjecture.

\newtheorem*{Conjecture}{Conjecture}
\begin{Conjecture}
Let $l$, $m$ and $r$ be three nonzero complex numbers with
$m=\unsurr-r$. The CGW algebra with parameters $l$ and $m$ over the
field $\Q(l,r)$ is semisimple except possibly if $r$ is a root of
unity or if
$$l\in\bigcup_{k=4}^n\lb
r^{2k-5},-r^{5-2k},-r^{4k-7},r^{7-4k},r^{7-2k},-r^{2k-7}\rb$$
\end{Conjecture}

\noin The values of Theorem $2$ point $(ii)$ that don't depend on
$n$, id est, $r^3,-\unsur{r^3},-r$ and $\unsur{r^3},-r^3$ are
obtained with $k=4$ and $k=5$ respectively. As for the values that
depend on the integer $n$, they are obtained with $k=n$.\\

The paper is organized as follows. In section $2$ we introduce the
diagrammatic version of the CGW algebra of type $D_n$. These
diagrams were introduced by the authors in \cite{BAL}. They show
that the CGW algebra of type $D_n$ is isomorphic to a tangle algebra
of type $D_n$. We next use this tangle algebra to construct the
representation of $\cgwn$ announced in Theorem $1$ of the
introduction. In $\S3$, we show that when the representation is
reducible, the action on a proper invariant subspace is a
$\mcalh(D_n)$-action. We then investigate the existence of
irreducible $\mcalh(D_n)$-modules of small dimensions inside the
Cohen-Wales space. To finish the proof of reducibility, we use
induction on the integer $n$. In $\S4$, we finish proving the
theorems of the introduction.
\section{Construction}

\subsection{The tangle algebra of type $D_n$ of
Cohen-Gijsbers-Wales} In \cite{BAL}, Arjeh M. Cohen, Di\'e A.H.
Gijsbers and David B. Wales build a diagram algebra that they show
to be isomorphic to the CGW algebra of type $D_n$. Here is how the
elements $g_k$'s and $e_k$'s are represented in their tangle algebra.\\
\epsfig{file=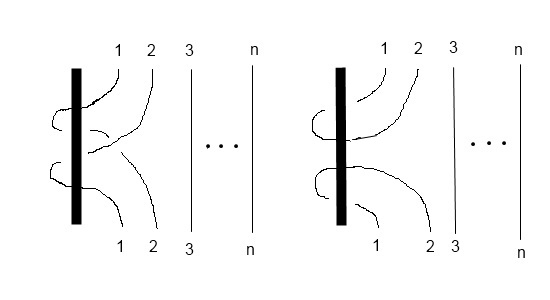, height=5.2cm}
\begin{center}\vspace{-1cm}\textit{The elements $g_1$ and
$e_1$}\end{center} \epsfig{file=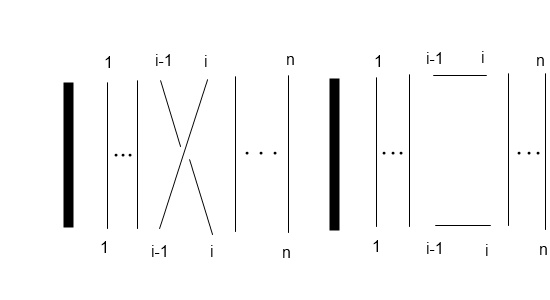, height=5.2cm}
\begin{center}\vspace{-1cm}\textit{The elements $g_i$ and $e_i$}\\
$2\leq i\leq n$\end{center} The vertical bar on the left hand side
is rigid and is called the pole. Among the $g_i$'s and the $e_i$'s,
only the elements $g_1$ and $e_1$ have strands twisting around the
pole. While the $g_i$'s only have vertical strands, the $e_i$'s
contain two horizontal strands. The following relations hold for
twists around the pole.\\ $\qquad\qquad$ \hspace{1cm}
\epsfig{file=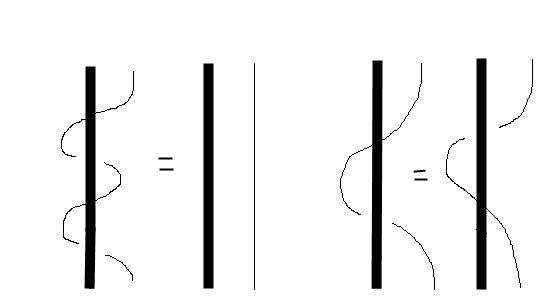, height=4.3cm}
\\The first relation says the pole has order two. The second
relation says it is indifferent whether the double twist is in the
order over-under or in the order under-over.
\newcounter{chiromains}\setcounter{chiromains}{2} The tangles are defined up to regular
isotopy, that is Reidemeister moves \Roman{chiromains} and
\setcounter{chiromains}{3} \Roman{chiromains} are permitted. For the
definition of these moves, see \cite{MOR}, page $4$. There are seven
more defining relations in the tangle algebra. Three of them are
independent of the presence of a first node and are the exact same
relations as in the algebra of Morton and Traczyk \cite{MOR}, the
diagrammatic version of the Birman-Murakami-Wenzl algebra. First,
there is the Kauffman skein relation. It is the way you transform an
under-crossing into an over-crossing and conversely. It is the
diagrammatic version of the algebraic equality
$g_i-g_i^{-1}=m\,(e_i-1)$.\\ \hspace{-3cm}\epsfig{file=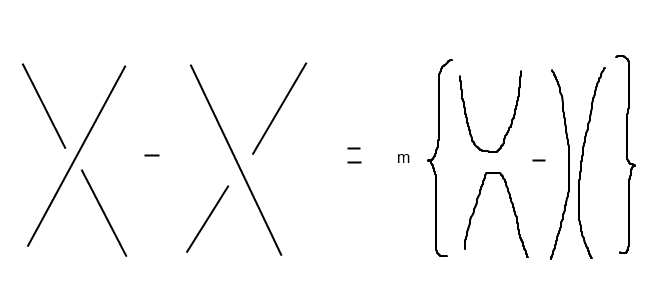,
height=4cm}\\ Next, the defining algebraic relations $(DL)$ are
called "delooping relations". Indeed, the way you get rid of a loop
on the diagrams is by multiplying by a factor $l$ or $l^{-1}$, as
follows.
\begin{center}\epsfig{file=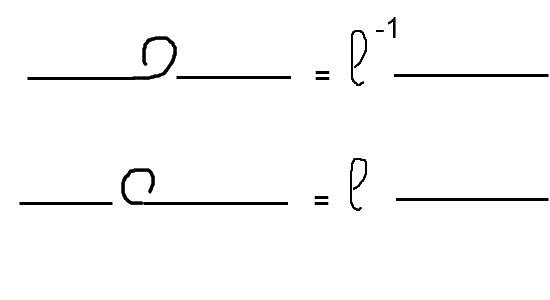,
height=4cm}\end{center}\vspace{-0.7cm} Finally, to finish with the
non pole-related relations, here is how the relation
$e_i^2=\delta\,e_i$ is conveyed in the diagrams. Each closed loop
not intersecting a tangle $T$ can be removed from the tangle by
multiplying by a factor $\delta$.
$$T\;\bigsqcup\;\bigcirc=\delta\;T$$
There are now four other relations that involve the pole. They are
the diagrammatic interpretations of the algebraic relations
$g_1\,g_2=g_2\,g_1$, $e_1\,g_2=g_2\,e_1$, $e_2\,g_1=g_1\,e_2$ and
$e_2\,e_1=e_1\,e_2$. We call the first of these relations the
commuting relation, as in \cite{BAL}. This relation will be
extensively used in the present
paper.\newcounter{chimaths}\setcounter{chimaths}{5} The other
relations are referred to by the authors in \cite{BAL} as the first
pole-related self-intersection relation, the second pole-related
self-intersection relation and the first closed pole loop relation
respectively. For these, we refer the reader to the diagrams
(\roman{chimaths}), \setcounter{chimaths}{6}(\roman{chimaths}) and
\setcounter{chimaths}{7}(\roman{chimaths}) of \cite{BAL}.
\begin{center}\epsfig{file=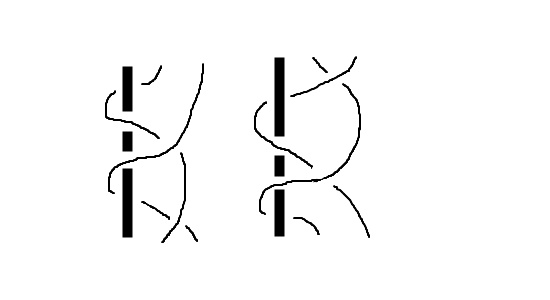,
height=4.8cm}\end{center}\vspace{-0.5cm}\textit{$\qquad\qquad\qquad\qquad\qquad\;$
The commuting relation}
\subsection{Construction of the representation}
To each vector $\ovl{w_{s,t}}$ of the Cohen-Wales space, we
associate a tangle. This tangle has two horizontal lines, one at the
top joining nodes $s$ and $t$ and one at the bottom joining nodes
$n-1$ and $n$. The top horizontal line over-crosses all the vertical
strands that it intersects. Moreover, if the vector wears a hat, the
top horizontal line twists around the pole, while when the vector
does not carry a hat, there is no twist around the pole. If there is
one twist around the pole, there should be another twist around the
pole. The first possible vertical strand twists around the pole with
the twist taking place below the twist of the horizontal strand.\\
In algebraic terms, to the root $\al_1$, one associates the CGW
algebra element $e_1\,e_{3,n}$. We then build the other positive
roots inductively by acting with the simple reflections, except for
the positive roots of type $w_{1,j}$. For instance,
$$\al_1+\al_2+2\,\al_3+\dots+2\,\al_i+\al_{i+1}+\dots+\al_j=r_i\dots
r_2\,r_j\dots r_3(\al_1)$$ and the associated agebra element is
$$g_{i,2}\,g_{j,3}e_1e_{3,n}$$
By $g_{i,j}$ or $e_{i,j}$, we understand $g_i\dots g_j$ or $e_i\dots
e_j$. The reader can check that the corresponding tangle has all the
above characteristics. For the positive roots of type $w_{1j}$, we
must also use some inverses of the $g_k$'s. For instance, when
$j\geq 3$, the associated CGW algebra element is
$g_2^{-1}g_1g_{j,3}e_1e_{3,n}$. An action by a generator $g_k$ on
these tangles can shift one of the horizontal strand's extremities
and/or introduce crossings between the vertical strands. Let
$\mcalh_n$ be the Hecke algebra of type $A_1\times D_{n-2}$ with
generators $z$ and $g_1,\,\dots,\,g_{n-2}$. Denote by $C_n$ the CGW
algebra of type $D_n$.
\newtheorem{Claim}{Claim}
\begin{Claim}
$M_n=C_n\,e_n/\langle C_n\,e_ie_j\,C_n\cap
C_n\,e_n\rangle_{i\not\sim j}$ is a right $\mcalh_n$-module for the
action: $$\begin{array}{l}
\forall\, 1\leq k\leq n-2,\forall\, x\in M_n,\, x\,.\,g_k=x\,g_k\\
\forall\, x\in
M_n,\,x\,.\,z=\unsur{\delta^2}\,\,x\,\,e_{n,3}\,e_1\,g_2\,e_1\,e_{3,n}
\end{array}$$
The $g_k$'s act to the right of elements in $M_n$ by simply
multiplying them to the right in $M_n$.\\
$z$ acts to the right of elements in $M_n$ by multiplying them to
the right by $\xi$ in $M_n$, where
$$\xi=\unsur{\delta^2}\,\,e_{n,3}\,e_1\,g_2\,e_1\,e_{3,n}$$
\end{Claim}
\noin\textsc{Proof of the Claim.} If $x\in M_n$, since for every
integer $k$ with $1\leq k\leq n-2$, the generator $g_k$ commutes to
$e_n$, we see that $x\,.\,g_k$ is again in $M_n$. Next, we have
$$x\,.\,(g_k^2+m\,g_k)=x\big(1+\frac{m}{l}\,e_k\big)$$
But since $e_n$ and $e_k$ commute and
$e_k\,e_n=0\;\text{in}\;C_n\,e_n/\langle C_n\,e_i\,e_j\,C_n\cap
C_n\,e_n\rangle_{i\not\sim j}$, we see that $g_k^2+m\,g_k$ acts like
the identity on $x$. Further, $x\,.\,z$ belongs to $M_n$. It remains
to show that $z^2+m\,z$ acts like the identity on $x$. We have
\begin{eqnarray}
\nts\nts\nts\nts x\,.\,(z^2+m\,z)&=&x\big(\unsur{\delta^4}\,e_{n,3}e_1g_2e_1e_{3,n}e_{n,3}e_1g_2e_1e_{3,n}+\frac{m}{\delta^2}\,e_{n,3}e_1g_2e_1e_{3,n}\big)\\
&=&x\big(\unsur{\delta^2}\,e_{n,3}\,e_1(g_2^2+m\,g_2)e_1e_{3,n}\big)\\
&=&x\big(\unsur{\delta^2}\,e_{n,3}\,e_1(1+\frac{l}{m}e_2)e_1e_{3,n}\big)\\
&=&x\big(\unsur{\delta}\,e_n\big)
\end{eqnarray}
Equality $(17)$ comes from the definition of the action. Equality
$(18)$ can be obtained by first using the relation
$e_n^2=\delta\,e_n$, then applying the relation $e_ie_je_i=e_i$ for
adjacent nodes $i$ and $j$ multiple times, finally by using the fact
that $e_1$ and $g_2$ commute and applying $e_1^2=\delta\,e_1$. To
get $(20)$, observe that $e_1e_2=0$ in $M_n$, then apply the same
machinery as before. Now $\unsur{\delta}\,e_n$ acts to the right
like the identity on any word ending in $e_n$. This settles the
claim.\hfill$\square$\\
Let's now provide our ground field $F=\Q(l,r)$ with a structure of
left $\mcalh_n$-module. We will consider the one-dimensional action
given by $g_k\,.\,1=r$ for every integer $k$ with $1\leq k\leq n-2$
and by $z\,.\,1=r$. Then,
$$C_n\,e_n/\langle C_n\,e_i\,e_j\,C_n\cap
C_n\,e_n\rangle_{i\not\sim
j}\otimes_{\mcalh_n}\Q(l,r)\in\,_{C_n}\text{Mod}$$ Our
representation is built inside this CGW algebra left module. We show
that the elementary tensors $\ovl{w_{s,t}}\otimes_{\mcalh_n} 1$ are
invariant under the action by the generators $g_k$'s. To do so, it
will be useful to understand the important role played by the
special element $\xi$. Here is how this element $\xi$ is represented
in the tangle algebra:\begin{center} \epsfig{file=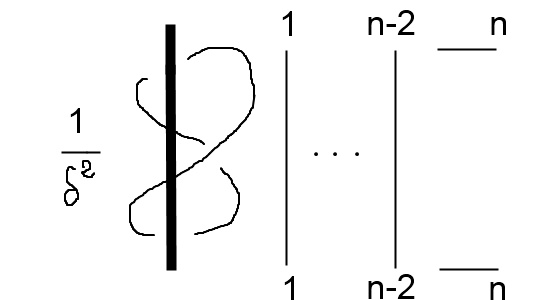,
height=3cm}\end{center} The strand which has the shape of an eight
can freely slide along the pole and it can be viewed as a
coefficient. It is called $\Xi^{+}$ in \cite{BAL} and it has many
interesting properties. One of them is that it commutes with another
twist around the pole, as shown on Fig. $9$ of \cite{BAL}. Another
property of $\Xi^{+}$ is that it satisfies a Kauffman skein type
relation (see equation $(2.1)$ of Lemma $2.11$ in \cite{BAL}). We
have $\xi=\unsur{\delta^2}\,\Xi^{+}\,e_n$. Note this simplified
expression for $\xi$ allows us to recover the fact from earlier that
$$\xi^2+m\,\xi=\unsur{\delta}\,e_n$$
Indeed,
\begin{eqnarray}
\xi^2+m\,\xi&=&\unsur{\delta^3}\,(\Xi^{+})^2\,e_n+\frac{m}{\delta^2}\,\Xi^{+}\,e_n\\
&=&\unsur{\delta^3}(\delta^2\,e_n-m\,\delta\,\Xi^{+}\,e_n+\frac{m}{l}\,\delta\,\Theta\,e_n)+\frac{m}{\delta^2}\Xi^{+}\,e_n\\
&=&\unsur{\delta}\,e_n\end{eqnarray} Equality $(21)$ holds by
definition. Equality $(22)$ comes from an application of equality
$(2.2)$ of Lemma $2.11$ of \cite{BAL}. Cohen, Gijsbers and Wales
define a $(0,0)$-tangle $\Theta$ that consists of two separate loops
each of which twists around the pole. The elements $\Xi^{+}$ and
$\Theta$ are interesting in that there is an analogy between $g_i$
(any $i$) and $\Xi^{+}$ on one hand and $e_i$ (any $i$) and $\Theta$
on the other hand, as is visible on equalities $(2.1)$ and $(2.4)$
of the same lemma. To get equality $(23)$, notice that
\begin{equation}\Theta\,e_n=e_{n,3}\,e_1\,e_2\,e_1\,e_{3,n}\end{equation}
\noin This equality is illustrated on the following
figure.\begin{center} \epsfig{file=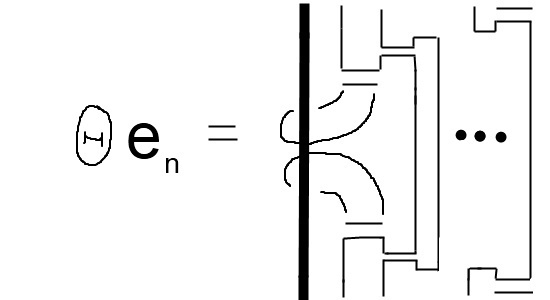,
height=4cm}\end{center} Now the right member of $(24)$ is zero in
$M_n$, which after simplification yields $(23)$. Let's now mention
the key property that we will use extensively to build the
representation. A loop around the pole can be suppressed at the cost
of a factor $\delta^{-1}\,\Xi^{+}$, as shown on Figure $1$.
\\\begin{center} \epsfig{file=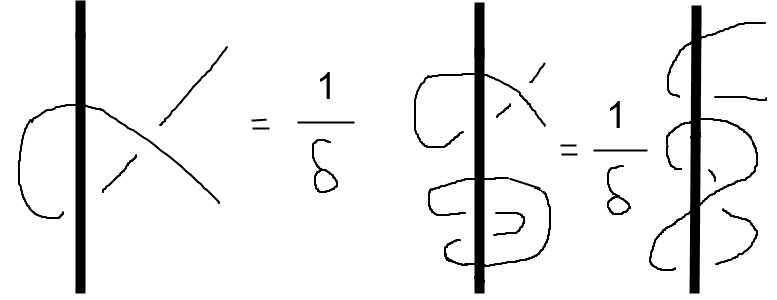, height=4cm}\end{center}\hfill \textit{Figure $1$}\\\\ The
trick is to add a closed loop by dividing by a factor $\delta$ and
to twist it twice around the pole, using the double twist relation.
Next, apply the first pole-related self-intersection relation
\setcounter{chimaths}{5}(\roman{chimaths}) of \cite{BAL} to get the
member to the right. In terms of our representation, here is how
this tangle property is nicely used. Replace $e_n$ by
$\unsur{\delta}\,e_n^2$. Then the CGW algebra element representing
the new tangle where the loop around the pole has been removed is
obtained from the old CGW algebra element representing the original
tangle containing a loop around the pole by multiplying it to the
right by $\unsur{\delta^2}\Xi^{+}\,e_n$. So suppressing the loop
around the pole is equivalent to acting to the right by $z$. But
recall we are working inside $M_n\otimes_{\mcalh_n}F$, so acting to
the right by $z$ on an element of $M_n$ is like acting to the left
by $z$ on $1$. So, in our representation, a loop around the pole is
replaced by a multiplication by $r$. We next show that, if in the
loop the crossing has the opposite sign, the loop can be removed at
the cost of a factor $\unsurr$. First by the same tangle trick, such
a loop can be removed at the cost of a factor
$\unsur{\delta}\;\Xi^{-}$ where $\Xi^{-}$ is the following tangle
\begin{center} \epsfig{file=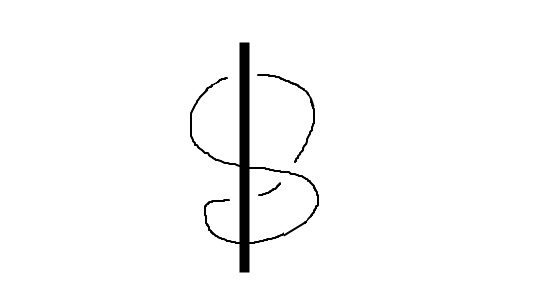,
height=4cm}\end{center}\begin{center}\textit{Figure $2$. The
$(0,0)$-tangle $\Xi^{-}$}\end{center} This tangle was introduced in
\cite{BAL}, and as mentioned by the authors, this tangle $\Xi^{-}$
is not the inverse of $\Xi^{+}$. We show that multiplying a tangle
by $\delta^{-1}\Xi^{-}$ is in the $C_n$-module
$M_n\otimes_{\mcalh_n}F$ a division by $r$. It suffices to show that
the product
$\big(\unsur{\delta^2}\,\Xi^{-}\,e_n\big)\big(\unsur{\delta^2}\,\Xi^{+}\,e_n\big)$
acts to the right of $e_n$ like the identity. We have
$$\bigg(\unsur{\delta^2}\,\Xi^{-}\,e_n\bigg)\bigg(\unsur{\delta^2}\,\Xi^{+}\,e_n\bigg)
=\unsur{\delta^3}\,\Xi^{-}\,\Xi^{+}\,e_n,$$ We show
$\Xi^{-}\Xi^{+}=\delta^2$. This is a simple application of the first
pole-related self-intersection relation, as shown on the figure
below.\begin{center}\epsfig{file=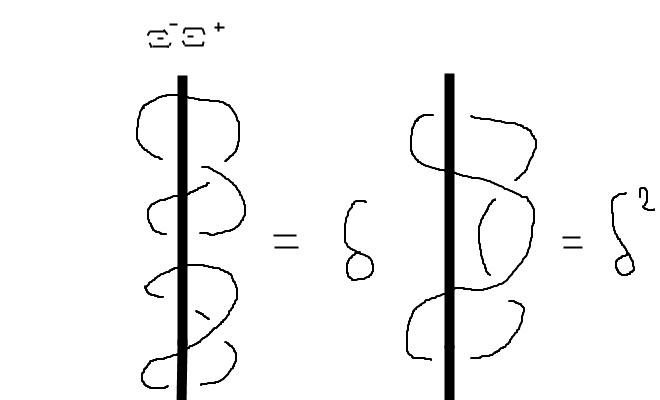,
height=4.8cm}\end{center}\vspace{.25cm} The "eight" at the bottom
becomes a closed loop, hence a factor $\delta$. The "eight" at the
top now has two self-intersections. These intersections are so that
an application of Reidemeister's move
\setcounter{chiromains}{2}\Roman{chiromains} is possible. Thus, you
get another closed loop and the announced result. So, the product we
considered is $\unsur{\delta}\,e_n$ and it indeed acts like the
identity to the right on $e_n$. We conclude that a loop around the
pole with a crossing of opposite sign as the one in Figure $1$ can
be removed at the cost of a division by $r$. \\
Using these preliminary remarks, it is now straightforward to see
that the action of the generators $g_k$'s leaves the basis
consisting of the elementary tensors
$\ovl{w_{s,t}}\otimes_{\mcalh_n}1$ invariant. The fact that we can
multiply at the bottom by $g_1$ (resp $g_1^{-1}$) at the cost of a
division by $r$ (resp a multiplication by $r$) allows us to easily
change the extremities of the vertical strand that twists around the
pole when these are not well positioned. This is for instance shown
on the following example. \begin{center} \epsfig{file=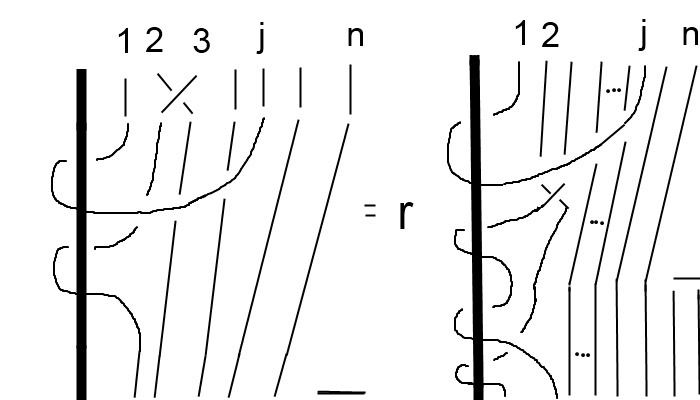,
height=4.8cm}\end{center}When computing the action by $g_3$ on the
basis vector $\wh{w_{1,j}}$, we use Reidemeister's move
\setcounter{chiromains}{3}\Roman{chiromains} to move the crossing
under the horizontal strand, then multiply the tangle at the bottom
by $g_1^{-1}$ and simultaneously compensate this addition by a
factor $r$. To finish, a simple use of the double twist relation,
followed by Reidemeister's move \setcounter{chiromains}{2}
\Roman{chiromains} allows node number $2$ at the top to join node
number $1$ at the bottom with a vertical strand twisting around the
pole below the horizontal strands. The final result is
$g_3\,.\,\wh{w_{1,j}}=r\,\wh{w_{1,j}}$.
We finish this construction section by describing one of the
actions, namely the action by $g_1$ on $\wh{w_{ij}}$ when $i\geq 3$.
This is the most complicated action for our representation.
Computing this action involves using the commuting relation a first
time, then the double twist relation once, Reidemeister's move
\setcounter{chiromains}{3}\Roman{chiromains} twice, then acting by
$g_2^{-1}$ at the bottom at the cost of a multiplication by $r$ to
get this tangle.
\\\epsfig{file=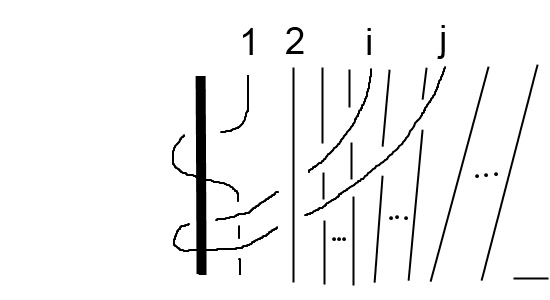, height=4.8cm}\\ The work
is not yet over. Indeed, in our basis, the horizontal strand always
twists above the vertical twist. An important feature of the
commuting relation is that it allows one to pull the bottom twist up
and to draw the upper twist down, hence changing the order in which
the horizontal strand and the vertical strand twist. In the process,
if the horizontal strand that twists around the pole was
over-crossing (resp under-crossing) the vertical strand that twists
around the pole, it now under-crosses (resp over-crosses) it. So,
using the commuting relation a second time, we
get\\
\epsfig{file=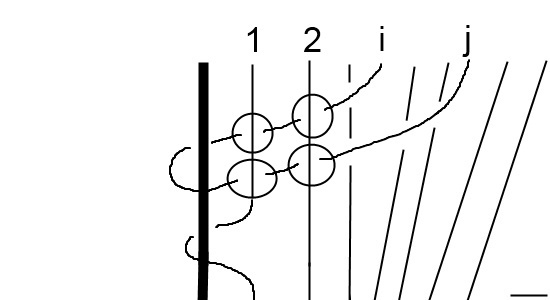, height=4.8cm}\\\\
We are now in a situation where the first two vertical strands
over-cross the top horizontal strand. We need to transform the four
crossings that are involved. These are pointed out on the diagram by
circles. We do so by using the Kauffman skein relation. When all the
under-crossings have been transformed into over-crossings, we get
the term $r\,\wh{w_{ij}}$ of $(1)$. The rest of the computations
must be done step by step with patience. They lead to the result in
$(1)$.\\
We last show another role played by the element $\xi$. When acting
by $g_1$ on $w_{12}$, one creates a pole-related self-intersection
as in Fig. $13$ of \cite{BAL}. As shown on the same figure, this
pole-related self-intersection can be replaced by
$\unsur{\delta}\,\Xi^{+}$. From there, replacing $e_n$ by
$\frac{e_n^2}{\delta}$, we get a multiplication to the right of
$e_n$ by $\xi$. Hence the action by $g_1$ on $w_{12}$ is a
multiplication by $r$.

\section{Reducibility}
\subsection{Action on a proper invariant subspace}
The following proposition is an easy but crucial statement about the
Cohen-Wales representation. It precises the action on a proper
invariant subspace when the representation is reducible.
\newtheorem{Proposition}{Proposition}
\begin{Proposition}
Let $\U$ be a proper invariant subspace of $V_n$. Then,
$$\forall\,1\leq i\leq n,\;\n^{(n)}(e_i)(\U)=0$$
Thus, the action on a proper invariant subspace is a Hecke algebra
action.
\end{Proposition}
\noin\textsc{Proof.} Let $\U$ be a proper invariant subspace of
$V_n$. Fix $i$ with $2\leq i\leq n$. An action by $e_i$ on any
element of the Cohen-Wales space is always proportional to the
vector $w_{i-1,i}$. Similarly, an action by $e_1$ on any element of
the Cohen-Wales space is always proportional to the vector
$\wh{w_{12}}$. Hence in the first case, if the action by $e_i$ on
$\U$ is non-trivial, then the vector $w_{i-1,i}$ belongs to $\U$.
And in the second case, if the action by $e_1$ on $\U$ is
non-trivial, then the vector $\wh{w_{12}}$ belongs to $\U$. But by
construction (see the beginning of $\S\,2.2$), the vectors
$\ovl{w_{s,t}}$'s are all of the form $y\,e_1e_{3,n}$ with $y$ a
certain product composed of $g_k$'s and $g_k^{-1}$'s. Then,
obviously if one of the $\ovl{w_{s,t}}$ belongs to $\U$, since $\U$
is invariant, $\U$ is then the whole space. This is in contradiction
with $\U$ proper. So, the proposition holds.\\\\
The goal now is to study which irreps of $\mcalh(D_n)$ can occur in
the Cohen-Wales space. First we recall some basic representation
theory of the Hecke algebra of type $D_n$ and study further the
degrees of the irreps of $\mcalh(D_n)$ that are less than $n^2-n$,
the degree of our representation $\n^{(n)}$.
\subsection{Degrees of the irreps of $\mcalh(D_n)$}
In this part, we assume that the Hecke algebra of type $D_n$ is
semisimple and we study the degrees of its irreducible
representations. We work over the field $F=\Q(l,r)$ which has
characteristic zero. Up to some rescaling of the generators, our
algebra $\mcalh(D_n)$ is the algebra $\mcalh_{r^2}(D_n)$ of
\cite{Hu1}. If $\mcalh(D_n)$ is semisimple, then by the proof of
Theorem $1.5$ of \cite{Hu1}, we have $f_n(r^2)\neq 0$. Pallikaros
defines $f_n(r^2)$ in his definition $2.12$ of \cite{PAL} as
$$f_n(r^2)=2\prod_{k=1}^{n-1}(1+r^{2k})$$
If $f_n(r^2)\neq 0$, then $r^{2k}\neq -1$ for every integer $k$ with
$1\leq k\leq n-1$. We will make this assumption in the remainder of
this paper. In the past fourty years, many authors have studied the
representation theory of the Hecke algebra of type $D_n$ \cite{Hu1},
\cite{PAL}, \cite{HOEF}, \cite{RAM} to only cite a few of them. It
seems the study finds its origin in the canadian Ph.D. thesis of
P.N. Hoefsmit \cite{HOEF}. Our work is based on the existing
theories that classify the irreducible $\mcalh(D_n)$-modules. The
main result that we use has been copied here from \cite{RAM}. We use
however the notations of our own paper.
\begin{thm}\textbf{[Hoefsmit], $1974$}
The modules $S^{(\al,\be)}$, where $(\al,\be)$ runs over all
unordered pairs of partitions such that $\al\neq\be$ and
$|\al|+|\be|=n$ and, when $n$ is even the modules
$S^{(\al,\al)^{+}}$ and $S^{(\al,\al)^{-}}$, where $\al$ runs over
all partitions such that $2|\al|=n$, form a complete set of
non-isomorphic irreducible modules for $\mcalh(D_n)$.
\end{thm}
\noin This Theorem says the non-isomorphic irreducible $\mcalh(D_n)$
are indexed by unordered double partitions of $n$. They are called
Specht modules. To keep the notations lighter, we will sometimes
write $S^{\al,\be}$ instead of $S^{(\al,\be)}$ and $S^{\al,\al^{+}}$
(resp $S^{\al,\al^{-}}$) instead of $S^{(\al,\al)^{+}}$ (resp
$S^{(\al,\al)^{-}}$), that is we omit the parenthesis around the
partitions. Since the double partitions are unordered, if
$(\la,\mu)$ is a double partition of $n$, we can assume without loss
of generality that $|\la|\leq|\mu|$ and so when $n$ is odd
$|\la|<\frac{n}{2}$. Moreover, the dimension of $S^{(\la,\mu)}$ is
given by the number of standard tableaux of shape $(\la,\mu)$,
except in the case when $n$ is even and $\la=\mu$ (and so
$2|\la|=n$). To describe the latter case, let's introduce the
notations of \cite{RAM}. When $\al\neq \be$, the Specht modules
$S^{(\al,\be)}$ are spanned by vectors $v_{L}$ indexed by standard
tableaux of shape $(\al,\be)$. When $\al=\be$, the Specht modules
$S^{(\al,\al)^{+}}$ and $S^{(\al,\al)^{-}}$ are respectively spanned
by vectors
$$\begin{array}{l} v_{L}^{+}=v_L+v_{\sigma\,L}\\
v_{L}^{-}=v_L-v_{\sigma\,L}
\end{array}$$
where $L$ is a standard tableau of shape $(\al,\al)$ and where
$\sigma$ is the map sending the standard tableau
$L=(L_{\al},L_{\be})$ of shape $(\al,\be)$ to the tableau
$\sigma\,L=(L_{\be},L_{\al})$ of shape $(\be,\al)$. Let's take an
example. Suppose $n=4$. The standard tableaux of shape $((2),(2))$
are: $$\begin{array}{l}L_1=\left(\;\begin{tabular}{|c|c|}\hline
1&2\\\hline\end{tabular}\;,\;\begin{tabular}{|c|c|}\hline
3&4\\\hline\end{tabular}\;\right)\qquad \sigma
L_1=\left(\;\begin{tabular}{|c|c|}\hline
3&4\\\hline\end{tabular}\;,\;\begin{tabular}{|c|c|}\hline
1&2\\\hline\end{tabular}\;\right)\\\\
L_2=\left(\;\begin{tabular}{|c|c|}\hline
1&3\\\hline\end{tabular}\;,\;\begin{tabular}{|c|c|}\hline
2&4\\\hline\end{tabular}\;\right)\qquad \sigma
L_2=\left(\;\begin{tabular}{|c|c|}\hline
2&4\\\hline\end{tabular}\;,\;\begin{tabular}{|c|c|}\hline
1&3\\\hline\end{tabular}\;\right)\\\\
L_3=\left(\;\begin{tabular}{|c|c|}\hline
1&4\\\hline\end{tabular}\;,\;\begin{tabular}{|c|c|}\hline
2&3\\\hline\end{tabular}\;\right)\qquad \sigma
L_3=\left(\;\begin{tabular}{|c|c|}\hline
2&3\\\hline\end{tabular}\;,\;\begin{tabular}{|c|c|}\hline
1&4\\\hline\end{tabular}\;\right)\end{array}$$

\noin We have
$\text{dim}(S^{(2),(2)^{+}})=\text{dim}(S^{(2),(2)^{-}})=3\;\text{and}\;
S^{(2),(2)}=S^{(2),(2)^{+}}\oplus S^{(2),(2)^{-}}$.\\
Recall the Cohen-Wales representation of $\cgwn$ has degree
$n(n-1)$, the number of positive roots of a root system of type
$D_n$. Our goal in this part is to find all the dimensions of the
Specht modules that have dimension less than $n(n-1)$ for a given
$n\geq 4$. We prove the following result.
\begin{theo}
Let $n$ be an integer with $n\geq 4$. Assume that $\mcalh(D_n)$ and
$\ih(n)$ are semisimple. \\
(i) Assume that $n\neq\lb 4,8\rb$. Then, the irreducible
representations of $\mcalh(D_n)$ have degrees
$$1,n-1,n,\cil,\frac{(n-1)(n-2)}{2}\;\;\text{or degrees greater
than}\;\;\frac{(n-1)(n-2)}{2}$$ (ii) The irreducible representations
of $\mcalh(D_4)$ have degrees $1,2,3,6$ or $8$.\\
(iii) The irreducible representations of $\mcalh(D_8)$ have degrees
$1,7,8,14,20,21$ or degrees greater than $21$. \\
(iv) For sufficiently large $n$, an irreducible representation of
$\mcalh(D_n)$ has degree\\\\ \indent
$1,n-1,n,\cil,\dbw,\chl,n(n-2)$\\\\ \indent or degree greater than
or equal to $n(n-1)$.
\end{theo}
\noindent\textsc{Proof.} Suppose $(\la,\mu)$ is a double partition
of $n$ with $|\la|=k$, $|\mu|=n-k$ and $n\geq 2k$. We study the
possible dimensions, depending on the value of $k$ and then $n$.
\begin{itemize}
\item[*] If $k=0$, $\la$ is the empty partition and $\mu$ is a
partition of $n$. We want to count the number of standard Young
tableaux of shape $\mu$. This number is the same as the dimension of
the Specht module $S^{\mu}$, where $S^{\mu}$ denotes a class of
irreducible $\ih(n)$-module for each partition $\mu$ of $n$. By
Corollary $2$ of \cite{CLDBW}, when $\ih(n)$ is semisimple, the
irreps of $\ih(n)$ have degree $1$, $n-1$, $\cil$, $\dbw$ or degrees
greater than $\dbw$, except when $n=4$ (resp $n=8$) when their
degrees are $1,2,3$ (resp $1,7,14,20,21$ or degree greater than
$21$). 
\item[*] If $k=1$, the Ferrers diagram of $\la$ is just one box.
There are $n$ possible choices to fill it. Once this single box is
filled, there are $(n-1)$ integers to fill in a standard tableau of
size $(n-1)$. Notice that $\frac{n(n-1)(n-4)}{2}\geq n(n-1)$ as soon
as $n\geq 6$. So, using the previous case, when $n\geq 6$ and $n\neq
9$, the only possible degrees are $n$ and $n(n-2)$. The case $n=9$
is in fact not an exception since $9\times 14>9\times 8$. When
$n=5$, the possible degrees are $5$, $10$ and $15$. Finally, when
$n=4$, the possible degrees are $4$ and $8$, so the case $n=4$ is
not exceptional.
\item[*] If $k=2$, suppose first $n\geq 5$. There are two possible partitions for $\la$ and
$\frac{n(n-1)}{2}$ ways to fill in the two boxes. Once this is
achieved, there are $(n-2)$ integers to fill in a standard tableau
of size $(n-2)$. Suppose first $n\neq 6$ and $n\neq 10$. When $n\geq
5$, we have $\frac{n(n-1)(n-3)}{2}\geq n(n-1)$, so the only
possibility is to have a degree equal to $\frac{n(n-1)}{2}$.
Consequently also, the case $n=10$ is not an exception. As when
$n=6$, because $2\,\frac{n(n-1)}{2}=n(n-1)$, the only possible
degree is again $15$. So this case is also not an exception. Let's
now deal with the case $n=4$. We have $|\la|=|\mu|=2$. Suppose first
$\la\neq\mu$. Since $(\la,\mu)$ is unordered, without loss of
generality, take $\la=(2)$ and $\mu=(1,1)$. The dimension of
$S^{(2),(1,1)}$ is $\binom{4}{2}=6$. When $\la=\mu$, the dimensions
of $S^{(2),(2)^{+}}$, $S^{(2),(2)^{-}}$ and their conjugates are $3$
as shown on the example above. Hence we have the extra degree $3$ in
that case.
\item[*] If $k=3$ and $n\geq 6$, we have $\binom{n}{3}=\frac{n(n-1)(n-2)}{6}\geq n(n-1)$
as soon as $n\geq 8$, hence we only need to worry about the cases
$n=6$ and $n=7$. First, if $n=6$ and $\la=\mu=(3)$, we have
$$\text{dim}(S^{(3),(3)^{+}})=\text{dim}(S^{(3),(3)^{-}})=\frac{\binom{6}{3}}{2}=10$$
Their conjugates also have dimension $10$ when $\la=\mu=(1,1,1)$. As
when $\la=\mu=(2,1)$, the dimension is too big. Indeed, we have
$$\text{dim}(S^{(2,1),(2,1)^{+}})=\text{dim}(S^{(2,1),(2,1)^{-}})=\frac{\binom{6}{3}\times
2\times 2}{2}=40>30$$
\\
Next, if $|\la|=|\mu|$ but $\la\neq\mu$, there are $\binom{3}{2}=3$
ways to choose two tableaux out of three. First if
$(\la,\mu)=(3),(1,1,1)$, the dimension of $S^{(\la,\mu)}$ is $20$.
Second if $\la=(2,1)$, the dimension is too big:
$2\times\binom{6}{3}=40>30=6\times 5$. \\
When $n=7$, we have $2\times\binom{7}{3}=70>42=7\times6$, so we must
have $\la\in\lb (3),(1^3)\rb$ and $\mu\in\lb (4),(1^4)\rb$. The
dimension of $S^{(\la,\mu)}$ is then $35$.
\item[*] If $k=4$ and $n\geq 8$, we have
$$\binom{n}{4}=\frac{n(n-1)(n-2)(n-3)}{6\times 4}\geq n(n-1)$$
Hence, except possibly when $n=8$, there are no irreducible
$\mcalh(D_n)$-modules of dimension less than $n(n-1)$. When $n=8$,
the Specht modules $S^{(4),(4)^{+}}$, $S^{(4),(4)^{-}}$ and their
two respective conjugates all have dimension $35$ and there are the
only irreducible $\mcalh(D_8)$-modules of dimension less than $56$.
\item[*] If $5\leq k\leq \frac{n}{2}$, the following inequality
holds.
\newtheorem{Lemma}{Lemma}
\begin{Lemma} For every integers $k$ and $n$ such that $n\geq 10$
and $5\leq k\leq \frac{n}{2}$, we have
$$\unsur{n(n-1)}\,\binom{n}{k}\geq 1$$
\end{Lemma}
\textsc{Proof of the Lemma.} The member to the left of the
inequality is
$$\frac{(n-2)\dots (n-k+1)}{k\dots 4\,.\,6}$$
In this fraction, there are $(k-2)$ terms in the nominator and there
are $(k-2)$ terms in the denominator. Moreover, we have
$n\geq\frac{n}{2}+2$, so that $n\geq k+2$. Further, we have
$n-k+1\geq 6$ as $n\geq \frac{n}{2}+5$.\hfill $\square$\\\\ When
$k<\frac{n}{2}$, the smallest possible dimension is $\binom{n}{k}$,
which by the lemma is greater than or equal to $n(n-1)$. Thus, it
remains to inspect the case when $k\geq 5$ and $n=2k$. In that case,
the smallest possible dimension is $\frac{1}{2}\binom{n}{k}$. But
when $n\geq 12$, we have $n\geq\frac{n}{2}+6$ so that
$n-\frac{n}{2}+2\geq 8$. Then,
$$\frac{1}{2}\frac{\binom{n}{\frac{n}{2}}}{n(n-1)}=\frac{(n-2)\dots
(n-\frac{n}{2}+2)(n-\frac{n}{2}+1)}{\frac{n}{2}\dots
5\,.\,8\,.\,6}\geq 1$$ and so
$$\frac{1}{2}\binom{n}{\frac{n}{2}}\geq n(n-1)$$
And when $n=10$, this inequality still holds by a direct
computation.\\ We conclude that when $5\leq k\leq \frac{n}{2}$,
there does not exist any irreducible $\mcalh(D_n)$-modules that have
dimension less than $n(n-1)$.
\end{itemize}
In summary, we have found the following degrees in the generic case:
$$1,n-1,n,\cil,\dbw, \chl,n(n-2)$$ If the study is complete for the non-zero $k$'s, it is incomplete when $k=0$.
Indeed, we don't have a complete list of the degrees of the irreps
of $\ih(n)$ when these degrees are between $\dbw$ and $n(n-1)$. We
only know that when $n$ is large enough, an irreducible
$\ih(n)$-module either belongs to $R_n(3)$ or has dimension greater
than $n^3$. This result comes from Theorem $5$ of \cite{JAMES},
applied with $m=3$. James' work deals with the irreducible
representations of the symmetric group $Sym(n)$, but can be applied
to the irreducible representations of $\ih(n)$. Indeed, in
characteristic zero when $\ih(n)$ is semisimple, the degrees of the
irreps of $Sym(n)$ are the same as the degrees of the irreps of
$\ih(n)$. James denotes by $R_n(m)$ the classes of irreducible
Specht modules $S^{\mu}$ with $\mu_1\geq n-m$, where
$\mu=(\mu_1,\mu_2,\dots)$ is a partition of $n$, or their
conjugates. A Specht module $S^{\mu}$ belongs to $R_n(3)$ if the
first row or the first column of the Ferrers diagram of the
partition $\mu$ contains $n-3$, $n-2$, $n-1$ or $n$ boxes. A
straightforward application of the Hook formula (see for instance
\cite{SAG}) shows that
$$\begin{array}{ccccc}
M\in R_n(1)&\Rightarrow&\text{dim}\,M\in&\nts\nts\lb 1,n-1\rb&(a)\\&&&\\
M\in R_n(2)\setminus
R_n(1)&\Rightarrow&\text{dim}\,M\in&\nts\nts\lb\cil,\dbw\rb&(b)\\&&&\\
M\in R_n(3)\setminus R_n(2)&\Rightarrow&\text{dim}\,M\in&\lb
\frac{n(n-1)(n-5)}{6},\frac{n(n-2)(n-4)}{3},\frac{(n-1)(n-2)(n-3)}{6}\rb&(c)\end{array}$$
For $(c)$, we give below the Ferrers diagrams and the hook
lengths. The first row of each diagram contains $(n-3)$ boxes.\\\\
\begin{tabular}[c]{|c|c|c|c|c|c|}
\hline n-2&n-3&n-4&n-6&\dots&1\\\hline 3&2&1\\\cline{1-3}
\end{tabular} $\;$
\begin{tabular}[c]{|c|c|c|c|c|}
\hline n-1&n-3&n-5&\dots&1\\\hline 3&1\\\cline{1-2} 1\\\cline{1-1}
\end{tabular} $\;$
\begin{tabular}[c]{|c|c|c|c|}
\hline n&n-4&\dots&1\\\hline 3\\\cline{1-1} 2\\\cline{1-1}
1\\\cline{1-1}
\end{tabular}\\\\
The dimensions of the respective Specht modules are obtained by
taking $n!$ over the product of the hook lengths. When $n\geq 11$,
the three quotients in $(c)$ are greater than or equal to $n(n-1)$.
Then, point $(iv)$ of Theorem $8$ holds. This settles the Theorem.
\subsection{Existence of a one-dimensional invariant subspace}
In this part we investigate the existence of a one-dimensional
invariant subspace inside $V_n$. We show the only values of $l$ and
$r$ for which that happens are those such that $l=\unsur{r^{4n-7}}$.
Assume such a space exists. Let $u$ be a spanning vector and so
there exists scalars $\la_1,\dots,\la_{n}$ such that
$\n_i(u)=\la_i\,u$ for each $i$. Further, since
$(\n_i^2+m\,\n_i)(u)=u$, we get $\la_i\in\lb r,-\unsurr\rb$. We show
that all the $\la_i$'s must in fact be equal to $r$. First, all the
$\la_i$'s are equal to say $\la$. Indeed, applying the braid
relation $g_1g_3g_1=g_3g_1g_3$, we get $\la_1=\la_3$ and applying
further the same braid relation on nodes $2,\dots,n$, we get that
all the $\la_i$'s with $2\leq i\leq n$ are equal. Let's write a
general form for $u$ as follows.
$$u=\sum_{1\leq i<j\leq n} \mu_{ij}\,w_{ij} +\sum_{1\leq i<j\leq
n}\wh{\mu_{ij}}\,\wh{w_{ij}}$$
\begin{Lemma}
Let $i$ be an integer with $2\leq i\leq n$. If $\n_i(u)=\la\,u$,
then $\ovl{\mu_{i,j}}=\la\,\ovl{\mu_{i-1,j}}$ for every $j>i$ and
$\ovl{\mu_{k,i}}=\la\,\ovl{\mu_{k,i-1}}$ for every $k<i-1$.
\end{Lemma}
\noindent \textsc{Proof.} Let's for instance prove the first
equality. It suffices to look at the term in $w_{i-1,j}$ (resp
$\wh{w_{i-1,j}}$) in $\n_i(u)=\la\,u$. For all $j>i$, an action by
$g_i$ creates a term in $w_{i-1,j}$ (resp $\wh{w_{i-1,j}}$) only
when it acts on $w_{i,j}$ (resp $\wh{w_{i,j}}$). Hence the
result.\\
\\As a corollary, if one of the $\mu_{ij}$'s (resp $\wh{\mu_{ij}}$'s) is
zero, then all the $\mu_{ij}$'s (resp $\wh{\mu_{ij}}$'s) are zero.
We show that it is impossible to have all the $\wh{\mu_{ij}}$'s
equal to zero. Indeed, if so, then all the $\mu_{ij}$'s are
non-zero. Acting with $\n_1$ on $w_{1,j}$ creates a term in
$\wh{w_{2j}}$. Moreover, as shown by the equations $(1)-(7)$, this
is the only way to create a term in $\wh{w_{2j}}$ when acting with
$\n_1$. This yields a contradiction. Thus, all the $\wh{\mu_{ij}}$'s
are non-zero. From there, it is easy to see that $\la$ must be equal
to $r$, not $-\unsurr$. Indeed, look at an action of $g_1$ on
$\wh{w_{34}}$ and notice this is the only way to get a term in
$\wh{w_{34}}$ when acting with $\n_1$. Since the term $\wh{w_{34}}$
is multiplied by $r$, we see that $\la$ must be equal to $r$. The
goal next is to find the relationship between the hat coefficients
and the non-hat coefficients. For that, we look at the coefficient
of $w_{1,j}$ for $j\geq 3$ in $\n_1(u)=r\,u$. We get
$$r\,\mu_{1j}=\wh{\mu_{2,j}}-m\,\mu_{1,j}-m\,\sum_{i=3}^{j-1}r^{i-3}\,\wh{\mu_{i,j}}-m\,\sum_{k=j+1}^{n}r^{k-5}\,\wh{\mu_{j,k}}$$
By simplifying this expression and using the relations between the
coefficients, we derive
$$\unsurr\,\mu_{1,j}=r\,\wh{\mu_{1,j}}-m\,\sum_{i=3}^{j-1}r^{i-3}r^{i-1}\,\wh{\mu_{1,j}}-m\,\sum_{k=j+1}^{n}r^{k-5}r^{k-1}\,\wh{\mu_{1,j}}$$
After evaluating the two sums of powers of $r$ and simplifying, we
obtain
\begin{equation}\mu_{1,j}=r^{2n-4}\,\wh{\mu_{1,j}}\end{equation} Let's
now look at the coefficient of $\wh{w_{12}}$ in $\n_1(u)=r\,u$. We
have
\begin{equation}\bigg(\unsur{l}-r\bigg)\,\wh{\mu_{12}}+m\,\sum_{j=3}^{n}r^{j-3}\,\mu_{2,j}+m\,\sum_{j=3}^{n}r^{j-4}\,\mu_{1,j}+m^2\bigg(\unsur{r^{11}}+\unsur{r^9}\bigg)\,\sum_{i=3}^{n-1}r^{2i}\sum_{j=i+1}^{n}r^{2j}\,\wh{\mu_{12}}=0\end{equation}
where we used the relation
$\wh{\mu_{ij}}=r^{i-1}r^{j-2}\,\wh{\mu_{12}}$. Also, we have
$$\begin{array}{l}\mu_{2j}=r^{2n-3}\,\wh{\mu_{1j}}=r^{2n-5+j}\,\wh{\mu_{12}}\\\\\mu_{1,j}=r^{2n-4}\,\wh{\mu_{1j}}=r^{2n-6+j}\,\wh{\mu_{12}}\end{array},$$
where the equalities to the left hold by $(25)$. After evaluating
the sums, simplifying and dividing by $\wh{\mu_{12}}$ which is known
to be a non-zero scalar, all the terms in $(26)$ simplify nicely. It
yields
$$l=\unsur{r^{4n-7}}$$
Conversely, suppose $l=\unsur{r^{4n-7}}$ and let
$$u=\sum_{1\leq i<j\leq n}r^{i+j}\,\wh{w_{ij}}+r^{2n-4}\sum_{1\leq i<j\leq n}r^{i+j}w_{ij}$$
It is a tedious but straightforward verification that the $g_k$'s
act on $u$ by multiplying it by $r$. Theorem $3$ is thus proven.
\subsection{Existence of an irreducible $(n-1)$-dimensional
invariant subspace} The goal of this part is to prove Theorem $4$
announced in the introduction. We still assume that $\mcalh(D_n)$
and $\ih(n)$ are semisimple. By $\S\,3.2$, except when $n=4$ and
when $n=6$, there are exactly two inequivalent irreducible
representations of $\mcalh(D_n)$ of degrees $(n-1)$. In \cite{THE},
page $53$, we provide matrix representations $(M_i)_{1\leq i\leq
n-1}$ for $S^{(n-1,1)}$ (resp $(N_i)_{1\leq i\leq n-1}$ for
$S^{(2,1^{n-2})}$) when we work with $\ih(n)$. To get a matrix
representation $(H_i)_{1\leq i\leq n}$ for $S^{(0),(n-1,1)}$ (resp
$(K_i)_{1\leq i\leq n} $ for $S^{(0),(2,1^{n-2})}$), it suffices to
take $H_{i+1}=M_i$ for all $i$ with $1\leq i\leq n-1$ and $H_1=H_2$
(resp $K_{i+1}=N_i$ for all $i$ with $1\leq i\leq n-1$ and
$K_1=K_2$). We show that it is impossible to have a basis
$v_1,\dots,v_{n-1}$ of vectors of the Cohen-Wales space such that
$$\nts\nts\nts\nts\nts\nts\nts\nts\nts\nts
\begin{array}{cc}(\nabla)&\begin{array}{cc}
\begin{array}{l}
(a)\;\n_1(v_1)=r\,v_1\\
(b)\;\n_1(v_2)=-r\,v_1-\unsurr\,v_2\\
(c)\;\n_1(v_t)=-\unsurr\,v_t\;\;\qquad\forall\,t\geq 3
\end{array}
&
\begin{array}{cc}
i\geq 2& \left|\begin{array}{ccccc}
(d)&\n_i(v_i)&=&-r\,v_{i-1}-\unsurr\,v_i&\\
(e)&\n_i(v_{i-1})&=&r\,v_{i-1}&\\
(f)&\n_i(v_{i-2})&=&-\unsurr\,(v_{i-2}+v_{i-1})&\\
(g)&\n_i(v_t)&=&-\unsurr\,v_t&\forall\,t\not\in\lb i,i-1,i-2\rb
\end{array}\right.\end{array}\end{array}\end{array}$$

\noindent In other words, the Specht module $S^{(0),(2,1^{n-2})}$
cannot occur in the Cohen-Wales space. First, we claim that for
$n\geq 8$, the result is obvious. Indeed, we have by equation $(g)$
of $(\nabla)$ \begin{equation}\forall t\geq
4,\,\n_t(v_1)=-\unsurr\,v_1\end{equation} Then, for every $t\geq 4$,
all the terms in $v_1$ must have indices starting or ending in $t-1$
or $t$. This is not possible as soon as $n\geq 8$. Thus, it remains
to deal with the cases $n\in\lb 4,5,6,7\rb$. When $n=7$, the
contradiction comes almost immediately as by $(27)$, we must have
$v_1=\la_{46}\,w_{46}+\wh{\la_{46}}\,\wh{w_{46}}$. But then
$\n_7(v_1)\neq -\unsurr\,v_1$. When $n=6$, we have by the same
arguments as above
$$v_1=\ovl{\la_{35}}\,\ovl{w_{35}}+\la_{45}\,w_{45}+\ovl{\la_{46}}\,\ovl{w_{46}}$$
By $\n_4(v_1)=-\unsurr\,v_1$, we can reduce further the expression
to the first three terms:
$$v_1=\ovl{\la_{35}}\,\ovl{w_{35}}+\la_{45}\,w_{45}$$ Then acting with $\n_6$ closes the case.
When $n=5$, by $(27)$ with $t\in\lb 4,5\rb$, we have
\begin{equation}
v_1=\la_{45}\,w_{45}+\la_{34}\,w_{34}+\ovl{\la_{35}}\,\ovl{w_{35}}
\end{equation}
Apply Lemma $2$ with $\la=-\unsurr$ and $i=5$ to get
$\wh{\la_{35}}=0$ and $\la_{35}=-\unsurr\,\la_{34}$. Apply again
Lemma $2$ with $\la=-\unsurr$ and $i=4$ to further get
$\la_{45}=-\unsurr\,\la_{35}$. Because the three coefficients that
are involved are thus related, they are all non-zero. By the
linearity in the $v_i$'s in the relations $(\nabla)$, we can set
without loss of generality $\la_{34}=1$. Then,
\begin{equation}
v_1=w_{34}-\unsurr\,w_{35}+\unsur{r^2}\,w_{45}
\end{equation}
We can now conclude. By $v_2=-r\,\n_3(v_1)-v_1$, there is no term in
$w_{15}$ in $v_2$. But there is a non-zero term in $w_{15}$ in
$\n_2(v_2)$, namely $w_{15}$. This contradicts
$\n_2\,v_2=-r\,v_1-\unsurr\,v_2$ and finishes the case $n=5$. Let's
deal with the case $n=4$. First, by $(27)$ with $t=4$, there are no
terms in $\wh{w_{34}}$ and $\ovl{w_{12}}$ in $v_1$. Hence, a general
form for $v_1$ is
\begin{equation}
v_1=\la_{34}\,w_{34}+\ovl{\la_{13}}\,\ovl{w_{13}}+\ovl{\la_{14}}\ovl{w_{14}}+\ovl{\la_{23}}\,\ovl{w_{23}}+
\ovl{\la_{24}}\,\ovl{w_{24}}
\end{equation}
We next apply Lemma $2$ with the relations $\n_4(v_1)=-\unsurr\,v_1$
and $\n_2(v_1)=r\,v_1$ to get the set of relations
$$\begin{array}{ccc}
\ovl{\la_{23}}&=&r\,\ovl{\la_{13}}\\
\ovl{\la_{14}}&=&-\unsurr\,\ovl{\la_{13}}\\
\ovl{\la_{24}}&=&-\ovl{\la_{13}}
\end{array}$$
Consequently also, at least one of $\la_{13}$ or $\wh{\la_{13}}$
must be non-zero. Otherwise, $v_1$ would be a multiple of $w_{34}$.
Then $v_2=-r\,\la(w_{24}+m\,w_{23}+r\,w_{34})$ for some non-zero
scalar $\la$. This expression is not compatible with
$\n_3(v_2)=r\,v_2$ since the term in $w_{23}$ in $\n_3(v_2)$ is
$-mr^2\la-\frac{m\,r\la}{l}$ and $\frac{m\,r\la}{l}\neq 0$. Further,
we see on the defining relations for the representation that $(5)$
is the only way to get a term in $\wh{w_{24}}$ when acting with
$\n_1$ on $v_1$. This fact together with $(a)$ implies that
$\la_{14}=r\,\wh{\la_{24}}$. Thus, both $\la_{13}$ and
$\wh{\la_{13}}$ are non-zero. Without loss of generality, set
$\wh{\la_{13}}=1$. Then, by $$v_2=-r\,\n_3(v_1)-v_1,$$ we see that
the coefficient of $\wh{w_{12}}$ in $v_2$ is $-r$. Thus, the
coefficient of $\wh{w_{12}}$ in $-r\,v_1-\unsurr\,v_2$ is $1$, while
the coefficient of $\wh{w_{12}}$ in $\n_2(v_2)$ is $-r^2$. Since
$r^2\neq -1$, this contradicts equality $(d)$ with $i=2$. So, we are
done with all the cases and conclude that the Specht module
$S^{(0),(2,1^{n-2})}$ cannot occur in the Cohen-Wales space. We now
show that the conjugate Specht module $S^{(0),(n-1,1)}$ can occur in
the Cohen-Wales space for the values $l=\unsur{r^{2n-7}}$. If in the
Cohen-Wales space there exists an irreducible invariant subspace
isomorphic to $S^{(0),(n-1,1)}$, there must exist a basis
$v_1,\dots,v_{n-1}$ such that the $v_i$'s satisfy the relations

$$\nts\nts\nts\nts\nts\nts\nts\nts\nts\nts
\begin{array}{cc}(\Delta)&\begin{array}{cc}
\begin{array}{l}
(a^{'})\;\n_1(v_1)=-\unsurr\,v_1\\
(b^{'})\;\n_1(v_2)=\unsurr\,v_1+\,v_2\\
(c^{'})\;\n_1(v_t)=r\,v_t\;\;\qquad\forall\,t\geq 3
\end{array}
&
\begin{array}{cc}
i\geq 2& \left|\begin{array}{ccccc}
(d^{'})&\n_i(v_i)&=&\unsurr\,v_{i-1}+r\,v_i&\\
(e^{'})&\n_i(v_{i-1})&=&-\unsurr\,v_{i-1}&\\
(f^{'})&\n_i(v_{i-2})&=&r\,(v_{i-2}+v_{i-1})&\\
(g^{'})&\n_i(v_t)&=&r\,v_t&\forall\,t\not\in\lb i,i-1,i-2\rb
\end{array}\right.\end{array}\end{array}\end{array}$$

\noindent The relations $(\Delta)$ are the conjugate relations of
$(\nabla)$ where $r$ has been replaced by $-\unsurr$. We show that
these relations force $l=\unsur{r^{2n-7}}$. Let's first use
$(e^{'})$ with $i=2$ to see that in $v_1$, there are no terms in
$\ovl{w_{s,t}}$ with $s\geq 3$ and there is no term in
$\wh{w_{12}}$. So,
\begin{equation}
v_1=\sum_{t=2}^n\mu_{1t}\,w_{1t}+\sum_{t=3}^n\wh{\mu_{1t}}\,\wh{w_{1t}}+\sum_{t=3}^n\mu_{2t}\,w_{2t}+\sum_{t=3}^n\wh{\mu_{2t}}\,\wh{w_{2t}}
\end{equation}
To get the explicit expression for $v_1$, it suffices now to juggle
with equations $(2)-(5)$. Look at the term in $w_{12}$ in
$\n_1(v_1)=-\unsurr\,v_1$ and get
$$r\,\mu_{12}+m\,\sum_{j=3}^n r^{j-3}\wh{\mu_{2j}}-m\,\sum_{j=3}^n r^{j-4}\,\mu_{1j}=-\unsurr\,\mu_{12}$$
Using $\mu_{1j}=-\unsurr\,\wh{\mu_{2j}}$, now derive
\begin{equation}
\mu_{12}=-\frac{m}{r}\sum_{t=3}^nr^{t-3}\,\wh{\mu_{2t}}
\end{equation}
For $t\geq 4$, we have $\n_t(v_1)=r\,v_1$. By Lemma $2$, it follows
that $\wh{\mu_{2t}}=r\,\wh{\mu_{2,t-1}}$ for every $t\geq 4$. Using
these relations in $(32)$ yields
\begin{equation}
\mu_{12}=\frac{\wh{\mu_{23}}}{r^2}(r^{2n-4}-1)
\end{equation}
Look at the coefficient of $\wh{w_{13}}$ in
$\n_1(v_1)=-\unsurr\,v_1$ and get using defining relations $(4)$ and
$(5)$
$$\mu_{23}=m\,\mu_{13}-\unsurr\,\wh{\mu_{13}}$$
Replace $\mu_{13}=-\unsurr\,\wh{\mu_{23}}$ and
$\wh{\mu_{13}}=-r\,\wh{\mu_{23}}$ (by $(e^{'})$ with $i=2$ and Lemma
$2$) to get \begin{equation}\mu_{23}=\unsur{r^2}\,\wh{\mu_{23}}
\end{equation}
All the coefficients in $v_1$ are now determined by $\wh{\mu_{23}}$.
Setting $\wh{\mu_{23}}=1$, we get the expression of $v_1$ given in
Theorem $4$:
\begin{equation}
v_1=(r^{2n-6}-\unsur{r^2})\,w_{12}+\sum_{j=3}^nr^{j-5}\,\big((w_{2j}-r\,w_{1j})+r^2(\wh{w_{2j}}-r\,\wh{w_{1j}})\big)
\end{equation}
Once $v_1$ is known, all the $v_i$'s for $2\leq i\leq n-1$ are
recursively determined by the formula $(f^{'})$. At this point, it
is tempting to look at the coefficient of $w_{12}$ in
$\n_2(v_1)=-\unsurr\,v_1$. However, the terms in $\unsur{l}$
simplify and this relation appears to be a
tautology. Thus, we cannot bypass involving the basis vector $v_2$. 
In fact, we show that the relations involving $v_1$ and $v_2$ are
enough to force a relationship between the parameters $l$ and $r$.
It suffices to look at the term in $w_{12}$ in $(d^{'})$ with $i=2$.
We derive \begin{equation}l=\unsur{r^{2n-7}}\end{equation}
Using inductively $(f^{'})$ with $(35)$ and replacing $l$ by its
value, we obtain the formulas of Theorem $4$.\\ To finish proving
the necessary condition of Theorem $4$, we still need to deal with
the special case $n=6$. When $n=6$ there are two more inequivalent
irreducible representations of $\mcalh(D_6)$ of degrees $5$. The
corresponding double partitions are $((0),(3,3))$ and
$((0),(2,2,2))$. We show that it is impossible to have an
irreducible $5$-dimensional invariant subspace that is isomorphic to
$S^{(0),(3,3)}$ or to $S^{(0),(2,2,2)}$. First we show the following
lemma.
\begin{Lemma} Let $n$ be an integer with $n\geq 5$.\\
(i) If in $V_n$ there exists an irreducible invariant subspace that
is isomorphic to the Specht module $S^{(0),(n-2,2)}$, then $l=r^3$.
When $n=5$, it is unique and it is spanned over $\Q(l,r)$ by the
vectors
\begin{equation}
v_4=r\,(w_{14}+r^2\,\wh{w_{14}})-(w_{24}+r^2\,\wh{w_{24}})+r(w_{23}+r^2\,\wh{w_{23}})-r^2\,(w_{13}+r^2\,\wh{w_{13}})\end{equation}\\
\vspace{-1cm}\begin{eqnarray}v_5&=&g_5\,.\,v_4\\
v_1&=&g_3\,.\,v_4-r\,v_4\\
v_2&=&g_3\,.\,v_5-r\,v_5\\
v_3&=&g_4\,.\,v_2-r\,v_2
\end{eqnarray}
(ii) In $V_n$, there does not exist any irreducible invariant
subspace that is isomorphic to the Specht module
$S^{(0),(2,2,1^{n-4})}$.
\end{Lemma}
\textsc{Proof of the Lemma.} In \cite{THE}, we found matrix
representations $(P_i)_{1\leq i\leq 4}$ and $(Q_i)_{1\leq i\leq 4}$
of degree $5$ for respectively $S^{(3,2)}$ and $S^{(2,2,1)}$. This
is Fact $1$ page $77$ of \cite{THE}. These are matrix
representations of $\ih(5)$. To get matrix representations of
$\mcalh(D_5)$, take $S_1=S_2=P_1$ and $S_i=P_{i-1}$ (resp
$T_1=T_2=Q_1$ and $T_i=Q_{i-1}$) for each $i\in\lb 3,4,5\rb$.

Let's prove $(i)$. We will need to use the branching rule. Branching
rules for Hecke algebras of type $D_n$ are stated in \cite{Hu2}.
Precisely, we use the results of Theorems $2.5$ and $2.6$ and
Corollary $2.8$. Jun Hu studies the decompositions into irreducible
modules of the socle of the restriction of each irreducible
$\mcalh(D_n)$-representation to $\mcalh(D_{n-1})$, that is for every
irreducible $\mcalh(D_n)$-module $D$, he describes
$Soc(D\da_{\mcalh(D_{n-1})})$. When we assume that $\mcalh(D_n)$ is
semisimple, the socle of a $\mcalh(D_n)$-module is the module
itself. Suppose $\W$ is an irreducible invariant subspace of $V_n$
that is isomorphic to $S^{(0),(n-2,2)}$. Then, by the branching
rule, the restriction of $\W$ to $\mcalh(D_5)$ is isomorphic to a
direct sum of Specht modules with one of the summands being
$S^{(0),(3,2)}$. The latter Specht module is obtained by $(n-5)$
successive removals of one box on the first row of the Ferrers
diagram of the partition $(n-2,2)$. Then, there exists in $\W$ a
family of five linearly independent vectors $(v_i)_{1\leq i\leq 5}$
such that the action by the $g_k$'s with $1\leq k\leq 5$ on these
vectors is given by the matrices $S_i$'s that were introduced above.
We work with these matrices to derive the results of $(i)$. We read
on the matrices $S_1$, $S_2$ and $S_4$ that
\begin{equation}g_i(v_4)=-\unsurr\,v_4\;\;\forall i\in\lb
1,2,4\rb\end{equation} These relations simplify greatly the shape of
$v_4$ and we immediately have
\begin{equation}v_4=\ovl{\la_{13}}\,\ovl{w_{13}}+\ovl{\la_{14}}\,\ovl{w_{14}}+
\ovl{\la_{23}}\,\ovl{w_{23}}+\ovl{\la_{24}}\,\ovl{w_{24}}\end{equation}
By the same relations $(42)$ and using Lemma $2$, the
hat-coefficients are all related, and so are the non-hat
coefficients. Moreover, by looking at the coefficient of
$\wh{w_{24}}$ in $(42)$ with $i=1$, the non-hat coefficients are
related to the hat coefficient by
$\la_{14}=-\unsurr\,\wh{\la_{24}}$. In particular, all the
coefficients involved in $(43)$ are non-zero. Without loss of
generality, we set $\la_{14}=r$. Then, we have
$$v_4=r\,w_{14}-w_{24}+r\,w_{23}-r^2\,w_{13}+r^2(r\,\wh{w_{14}}-\wh{w_{24}}+r\,\wh{w_{23}}-r^2\,\wh{w_{13}})$$
Up to a reordering of the terms, this is formula $(37)$ in the
statement of Lemma $3$. Once $v_4$ is known, the vectors
$v_1,v_2,v_3$ and $v_5$ are then uniquely determined by the formulas
$(38)-(41)$ which follow after a glance at the matrices $S_i$'s. The
uniqueness part when $n=5$ is then established. Further, we show
that $l$ and $r$ must be related in a certain way. Combining $(39)$
and the relation $\n_1(v_1)=r\,v_1+v_4$, we get
\begin{equation}
\n_1\n_3(v_4)-r\,\n_1(v_4)=r\,\n_3(v_4)+(1-r^2)\,v_4
\end{equation}
Look at the term in $\wh{w_{12}}$ in this expression. First, we
compute
\begin{multline}\n_3(v_4)=r^2(w_{14}+r^2\wh{w_{14}})-(w_{34}+r^2\,\wh{w_{34}})+\frac{r^5}{l}\,w_{23}+r^4\,\wh{w_{23}}\\+m\,r^2(w_{13}+r^2\,\wh{w_{13}})-r^2(w_{12}+r^2\,\wh{w_{12}})\end{multline}
Then we use defining relations $(1)$, $(4)$, $(5)$, $(6)$ of the
representation to get
$$r^2\,m-r^2m^2\big(\unsurr+r\big)+\frac{r^5}{l}m+mr^2\frac{m}{r}-\frac{r^4}{l}-r\,\big\lb rm-mr+rm-r^2\,\frac{m}{r}\big\rb=-r^5$$
This expression simplifies to yield $l=r^3$. This finishes the proof
of point $(i)$. Let's prove $(ii)$. Suppose $\W$ is an irreducible
invariant subspace of $V_n$ that is isomorphic to
$S^{(0),(2,2,1^{n-4})}$. Then, by the branching rule,
$S^{(0),(2,2,1)}$ is a component of $\W\da_{\mcalh(D_5)}$. Hence
there exists linearly independent vectors $v_1,\dots,v_5\,$ so that
the actions by $g_1,\dots,g_5$ on these vectors is given by the
matrices $T_1,\dots,T_5$ that were introduced at the beginning of
the proof. We show the relations force $v_1=0$, hence a
contradiction and the result. Denote the coefficients of $v_1$ by
$\ovl{\la_{i,j}}$. From
$$\left\lb\begin{array}{ccc}
g_2.\,v_1&=&-\unsurr\,v_1+v_4\\
g_3.\,v_4&=&v_1-\unsurr\,v_4
\end{array}\right.$$
derive \begin{equation}
g_3g_2.\,v_1+\unsurr\,g_3.v_1+\unsurr\,g_2.v_1=(1-\unsur{r^2})\,v_1
\end{equation}

Notice $$r^2+1+\unsur{r^2}=0\iff r=\frac{\pm 1\pm i\sqrt{3}}{2}$$
\noindent This is impossible when $(r^2)^3\neq 1$. So $(46)$ implies

$$\ovl{\la_{i,j}}=0\;\;\text{for all $i\geq 4$ and all $j\geq 5$}$$

\noindent Now, use the same trick a second time with
$$\left\lb\begin{array}{ccc}
g_4.\,v_1&=&-\unsurr\,v_1+v_4\\
g_3.\,v_4&=&v_1-\unsurr\,v_4
\end{array}\right.$$
and derive $$\ovl{\la_{1,j}}=0\;\;\text{for all $j\geq 5$}$$

\noindent We will keep making $v_1$ lighter. For now,

$$v_1=\ovl{\la_{12}}\ovl{w_{12}}+r\ovl{\la_{12}}\ovl{w_{13}}+\ovl{\la_{14}}\ovl{w_{14}}
+\sum_{j=3}^n\ovl{\la_{2,j}}\ovl{w_{2,j}}+r\,\sum_{j=4}^n\ovl{\la_{2,j}}\ovl{w_{3,j}}$$

\noindent where we also used the relation $g_3.\,v_1=r\,v_1$
together with Lemma $2$. Notice on the matrices that
\begin{equation}g_2.\,v_1=g_4.\,v_1\end{equation}
It follows that
$$\ovl{\la_{2,j}}=0\;\;\text{for all $j\geq 5$}$$ \noindent Now
$v_1$ reduces to
$$v_1=\ovl{\la_{12}}\ovl{w_{12}}+r\ovl{\la_{12}}\ovl{w_{13}}+\ovl{\la_{14}}\ovl{w_{14}}+\ovl{\la_{23}}\,\ovl{w_{23}}+\ovl{\la_{24}}\,\ovl{w_{24}}+r\,\ovl{\la_{24}}\,\ovl{w_{34}}$$
By looking at the term in $\ovl{w_{13}}$ in $(47)$, we further get
\begin{equation}\ovl{\la_{23}}=\ovl{\la_{14}}\end{equation}
And by looking at the term in $\ovl{w_{14}}$ in $(47)$, we get
\begin{equation}\ovl{\la_{24}}=r\,\ovl{\la_{12}}-m\,\ovl{\la_{14}}\end{equation}
We show that $\ovl{\la_{12}}=0$. Below, we write the relations that
we use and the relations that they imply on the coefficients. We
write $\ovl{\la_{s,t}^{(i)}}$ for the coefficient of $\ovl{w_{s,t}}$
in $v_i$. So $\ovl{\la_{s,t}^{(1)}}$ is simply $\ovl{\la_{s,t}}$. We
have
$$(\mathcal{R})\left|\begin{array}{ccc}
v_2=g_5.\,v_1&\Longrightarrow&\ovl{\la_{12}^{(2)}}=r\,\ovl{\la_{12}^{(1)}}
\\&&\\v_3=g_4.\,v_2+\unsurr\,v_2&\Longrightarrow& \ovl{\la_{12}^{(3)}}=(1+r^2)\,\ovl{\la_{12}^{(1)}}\\&&\\v_4=g_4.\,v_1+\unsurr\,v_1&\Longrightarrow&
\ovl{\la_{12}^{(4)}}=\big(r+\unsurr\big)\,\ovl{\la_{12}^{(1)}}
\\&&\\v_5=g_5.\,v_4&\Longrightarrow&\ovl{\la_{12}^{(5)}}=(1+r^2)\,\ovl{\la_{12}^{(1)}}
\end{array}\right.$$
Next, we read on the third column of the matrix $T_5$ that
\begin{equation}
g_5\,.\,v_3=\unsurr\,v_1+v_2-\unsurr\,v_3-\unsur{r^2}\,v_4-\unsurr\,v_5
\end{equation}
We look at the coefficient of $\ovl{w_{12}}$ in $(50)$ and we use
$(\mathcal{R})$. We obtain
\begin{equation}
2(r+\unsurr)+r^3+\unsur{r^3}=0\;\;\text{or}\;\;\ovl{\la_{12}^{(1)}}=0
\end{equation}
Let's solve the equation in $r$. With $X=\unsurr+r$, the equation is
equivalent to $$X^3-X=0$$ From before, we know that $X=1$ or $X=-1$
are impossible, as $(r^2)^3$ would then be $1$. Also, $X=0$ is
impossible since it leads to $r^2=-1$, which is excluded. Then, we
are forced to have $\ovl{\la_{12}}=0$ where we forgot the index
$(1)$ to conform to the notations of the beginning. Now, plug back
$(48)$ and $(49)$ into the expression for $v_1$ and get the newer
and simpler expression
\begin{equation}
v_1=\ovl{\la_{23}}(\ovl{w_{14}}+\ovl{w_{23}})-m\,\ovl{\la_{23}}(\ovl{w_{24}}+r\,\ovl{w_{34}})
\end{equation}
This is enough to conclude. Indeed, by looking at the terms in
$\wh{w_{24}}$ and $\wh{w_{13}}$ in the relation
$g_1.\,v_1=g_4.\,v_1$, we get the respective equations
$$\left\lb\begin{array}{l}
\la_{23}-(1+m^2+m^2\,r^2)\,\wh{\la_{23}}=0\\
\la_{23}-(1+m^2)\,\wh{\la_{23}}=0
\end{array}\right.$$
Then, all the coefficients in $v_1$ are zero, a contradiction.
Therefore, there does not exist any irreducible invariant subspace
in $V_n$ that is isomorphic to the Specht module
$S^{(0),(2,2,1^{n-4})}$ and point (ii) of the lemma is proven. Let's
go back to the proof of Theorem $4$. Suppose $\W$ is an irreducible
$5$-dimensional subspace of $V_6$ that is isomorphic to
$S^{(0),(3,3)}$. Then, applying the branching rule yields
$$\W\da_{\mcalh(D_5)}\simeq S^{(0),(3,2)}$$
Then, $\W$ is spanned over $F$ by vectors $v_1,\dots,v_5$ given in
equations $(37)-(41)$ of point (i) of Lemma $3$. A quick inspection
at these vectors shows that node number $6$ does not appear in them.
However, when acting with $\n_6$ on $v_5$, one creates terms that
end in node number $6$. This is in contradiction with the fact that
$\W$ is spanned by the $v_i$'s, $1\leq i\leq 5$. We conclude that
there does not exist any irreducible invariant subspace of $V_6$
that is isomorphic to $S^{(0),(3,3)}$. Suppose that there exists a
$5$-dimensional irreducible invariant subspace $\W$ of $V_6$ that is
isomorphic to $S^{(0),(2,2,2)}$. Then,
$$\W\da_{\mcalh(D_5)}\simeq S^{(0),(2,2,1)}$$
There exists vectors $v_1,\dots,v_5$ of $\W$ such that the action by
$g_1,\dots,g_5$ on these vectors is given by the matrices $T_i$'s,
$1\leq i\leq 5$. We get the same contradiction as in the proof of
Lemma $3$ point $(ii)$. We conclude that there does not exist any
irreducible $5$-dimensional invariant subspace of $V_6$ that is
isomorphic to $S^{(0),(2,2,2)}$, and so out of the four irreducible
representations of $\mcalh(D_6)$ of degree 5, only one of them can
occur in the Cohen-Wales space $V_6$ and this when
$l=\unsur{r^5}$.\\
The necessary condition of Theorem $4$ is now entirely proven for
$n\geq 5$. Conversely, suppose $l$ and $r$ are related as in $(36)$
and define vectors $v_i$'s, $1\leq i\leq n-1$ as in Theorem $4$.
Clearly, these vectors are linearly independent and we can check
that they satisfy all the relations $(\Delta)$. Then, they span an
irreducible $(n-1)$-dimensional invariant subspace inside $V_n$.
This ends the proof of Theorem $4$ in the case when $n\geq 5$. In
the case when $n=4$, as seen in $\S3.2$, there are four more
non-isomorphic irreducible $\mcalh(D_4)$-modules of dimension $3$,
namely $S^{(2,2)^{+}}$ and $S^{(2,2)^{-}}$ and their conjugates. 

\subsection{Existence of an irreducible $n$-dimensional invariant
subspace} The object of this section is to prove Theorem $5$
announced in the introduction. Given $n\geq 4$, our study in $\S3.2$
shows that except when $n=5$, there are exactly two distinct classes
of irreducible $\mcalh(D_n)$-modules of dimension $n$, namely the
Specht module $S^{(1),(n-1)}$ and its conjugate $S^{(1),(1^{n-1})}$.
When $n=5$, there are exactly four distinct classes of irreducible
$\mcalh(D_5)$-modules of dimension $5$, namely the ones above and
the Specht modules $S^{(0),(3,2)}$ and $S^{(0),(2,2,1)}$. The latter
Specht modules have been studied in the previous section. We proved
in Lemma $3$ of that section that $S^{(0),(3,2)}$ may occur when
$l=r^3$, while $S^{(0),(2,2,1)}$ can never occur. We found a matrix
representation for $S^{(1),(3)}$ and the proof of Theorem $5$ will
rely on it. We give this representation in the following theorem.
\begin{theo}
The matrices
$$H_1=\begin{bmatrix}r&0&0&0\\
-r^2+\unsur{r^2}&\unsur{r^3}&-\unsur{r^2}-\unsur{r^4}&0 \\
-r^3+\unsurr&-r^2+\unsur{r^2}&r-\unsur{r^3}-\unsurr&0\\
1-r^2&\unsurr-r&-\unsur{r^2}&r  \end{bmatrix}, H_2=\begin{bmatrix}
r&0&0&0\\0&r&0&0\\0&0&r&0\\0&0&1&-\unsurr\end{bmatrix}$$
$$H_3=\begin{bmatrix}r&0&0&0\\0&r&0&0\\0&1&-\unsurr&1\\0&0&0&r\end{bmatrix},
H_4=\begin{bmatrix}r-\unsurr&1&-\unsurr&0\\1&0&1&0\\0&0&r&0\\0&0&0&r\end{bmatrix}\\\\$$
define an irreducible matrix representation of $\mcalh(D_4)$ of
degree $4$.
\end{theo}
\noindent\textsc{Proof.} It is easy to visualize that $H_2$, $H_3$
and $H_4$ satisfy the usual braid relations on nodes $2,3,4$ and
that $H_i^2+m\,H_i=I$ for each $i\in\lb 2,3,4\rb$, where $I$ denotes
the identity matrix of size $4$. Further, we check, for instance
with Maple, that $H_1^2+m\,H_1=I$, $H_1H_3H_1=H_3H_1H_3$ and that
$H_1$ commutes to both $H_2$ and $H_4$. Hence these matrices define
a representation of $\mcalh(D_4)$ of degree $4$. It remains to show
that this representation is irreducible. Suppose there exists a
one-dimensional invariant subspace spanned by
$u=(u_1,u_2,u_3,u_4)^{tr}$. We must have $H_i\,u=r\,u$ for all
$i\in\lb 1,2,3,4\rb$. Next, we used Maple to solve this system of
equations and got $u=0$. So there does not exist any one-dimensional
invariant subspace. If we can show that there does not exist any
irreducible $2$-dimensional invariant subspace as well, then we are
done by using the semisimplicity of $\mcalh(D_4)$. Up to
equivalence, there is a unique irreducible representation of
$\mcalh(D_4)$ of degree $2$ and it is defined by the matrices
$$J_1=J_2=J_4=\begin{pmatrix}-\unsurr&1\\0&r\end{pmatrix}\qquad\qquad
J_3=\begin{pmatrix}r&0\\1&-\unsurr\end{pmatrix}$$ So, there exists
two non-zero linearly independent vectors $v_1$ and $v_2$ of
$\mathbb{C}^4$
so that \\\\
$\forall i\in\lb 1,2,4\rb,\,\;\left\lb\begin{array}{ccc}
H_iv_1&=&-\unsurr\,v_1\;\;\;\;\;\;\;\;(\star)_i\\H_iv_2&=&v_1+r\,v_2\end{array}\right.\qquad\begin{array}{ccc}H_3v_1&=&r\,v_1+v_2\\
H_3v_2&=&-\unsurr\,v_2\end{array}$ \\\\
Relation $(\star)_i$ applied with $i=2,4$ suffices to force $v_1=0$
by using for instance the first two rows of $H_2$, the last two rows
of $H_4$ and the fact that $r^2\neq -1$. Thus, we get a
contradiction. This ends the proof of Theorem $9$.\\
Suppose there exists in $V_n$ an irreducible $n$-dimensional
invariant subspace $\W$ that is isomorphic to $S^{(1),(n-1)}$.
Applying the branching rule $(n-4)$ times yields
$$\W\da_{\mcalh(D_4)}\simeq (n-4)\,S^{(0),(4)}\oplus S^{(1),(3)}\qquad\qquad(\divideontimes)$$
So there must exist vectors $v_1,v_2,v_3,v_4$ in $\W$ so that the
left actions by $g_1$, $g_2$, $g_3$ and $g_4$ on these vectors are
given by the matrices $H_i$'s of Theorem $9$. The computations are
technical. We sketch them here and leave the details to the reader.
First we read on the matrix $H_3$ that $g_3v_3=-\unsurr\,v_3$.
Hence, a general form for $v_3$ is
$$v_3=\la_{23}\,w_{23}+\ovl{\la_{12}}\ovl{w_{12}}-\unsurr\,\ovl{\la_{12}}\ovl{w_{13}}+\sum_{j=4}^n\ovl{\la_{2j}}\ovl{w_{2j}}-\unsurr\,\sum_{j=4}^n\ovl{\la_{2j}}\ovl{w_{3j}}$$
We now use $v_4$ to get more relations between these coefficient.
First, we have $g_4v_4=r\,v_4$, and so
$\ovl{\la_{24}^{(4)}}=r\ovl{\la_{23}^{(4)}}$. Second, with
$v_4=g_2v_3-r\,v_3$, we obtain
\begin{eqnarray*}
\la_{24}^{(4)}&=&-\unsurr\,\la_{24}\\
\la_{23}^{(4)}&=&-\unsurr\,(\la_{12}+\la_{23})
\end{eqnarray*}
and \hspace{-0.5cm}\begin{eqnarray*}
\wh{\la_{24}^{(4)}}&=&-\unsurr\wh{\la_{24}}\\
\wh{\la_{23}^{(4)}}&=&-\unsurr\wh{\la_{12}}
\end{eqnarray*}
So we get $\la_{24}=r(\la_{12}+\la_{23})$ and
$\wh{\la_{24}}=r\,\wh{\la_{12}}$. \\
Further, use $g_1\,v_4=r\,v_4$ and $v_4=g_2\,v_3-r\,v_3$ to derive
$$g_1g_2\,v_3=r\,g_1v_3+r\,g_2v_3-r^2\,v_3$$ Look at the coefficient in
$\wh{w_{1j}},\,j\geq 4$ in this relation and get
$$\forall j\geq 4,\,\wh{\la_{2j}}=-\la_{2j}$$
In particular, doing $j=4$ and using the relations above, we get
\begin{equation}
\wh{\la_{12}}+\la_{12}=-\la_{23}
\end{equation}
Let's use the relations above to write
\begin{multline}
v_3=\ovl{\la_{12}}\ovl{w_{12}}-\unsurr\,\ovl{\la_{12}}\ovl{w_{13}}-(\la_{12}+\wh{\la_{12}})\,w_{23}+r\,\wh{\la_{12}}
(\wh{w_{24}}-w_{24})+\wh{\la_{12}}(w_{34}-\wh{w_{34}})
\\+\sum_{j=5}^n\la_{2j}\big(w_{2j}-\unsurr\,w_{3j}\big)-\sum_{j=5}^n\la_{2j}\big(\wh{w_{2j}}-\unsurr\,\wh{w_{3j}}\big)
\end{multline}
We will now find more relations between these coefficients and a
relation involving $l$. We read on the third column of $H_4$ that
\begin{equation}\big(\unsur{r^2}+\unsur{r^4}\big)g_4v_3=-\unsurr\,g_4\big(\unsur{r^2}+\unsur{r^4}\big)v_2+\big(\unsur{r^4}+\unsur{r^2}\big)v_2+\big(\unsurr+\unsur{r^3}\big)\,v_3,\end{equation}
where we multiplied both sides by $\unsur{r^2}+\unsur{r^4}$ and
where we used that $v_1=g_4v_2$. Further, we read on the third
column of $H_1$ that
\begin{equation}\big(\unsur{r^2}+\unsur{r^4}\big)v_2=-\unsur{r^2}g_2v_3+(r-\unsur{r^3})v_3-g_1v_3,\end{equation}
where we used that $v_4=g_2v_3-r\,v_3$. Plugging $(56)$ into $(55)$
and simplifying now yields the following equation in $v_3$.
\begin{equation}\big(\unsur{r^2}+1\big)g_4v_3=\unsur{r^3}g_4g_2v_3+\unsurr\,g_4g_1v_3-\unsur{r^2}g_2v_3-g_1v_3+(r+\unsurr)v_3\end{equation}
By looking at the coefficient in $w_{24}$ in $(57)$, we obtain the
relation
\begin{equation}
-\wh{\la_{12}}+\unsur{r^2}\la_{12}=m\,\sum_{j=5}^n r^{j-6}\,\la_{2j}
\end{equation}
We will use this expression of the sum on the right hand side to
derive a relation involving $\la_{12}$, $\wh{\la_{12}}$ and $l$. It
suffices to look at the term in $w_{23}$ in $g_3v_3=-\unsurr\,v_3$.
We get, where we used $(53)$ and $(54)$,
\begin{equation}
\la_{23}\big(\unsur{l}+\unsurr+\frac{m}{lr}\big)=m\,\sum_{j=5}^n\la_{2j}r^{j-6}\big(r-\unsur{l}\big)-m\,\wh{\la_{12}}\big(1-\unsur{lr}\big)
\end{equation}
Replacing $\la_{23}$ as in $(53)$ and replacing the sum of $(59)$ as
in $(58)$, we then obtain
\begin{equation}
lr\,\la_{12}=-\wh{\la_{12}}
\end{equation}
Note we can already conclude in the case $n=4$. Indeed, in that case
it follows from $(58)$ that $-\wh{\la_{12}}+\unsur{r^2}\la_{12}=0$
and so we get $l=-\unsur{r^3}$ by using $(60)$. To solve the general
case, we introduce a few notations.
\begin{Claim}\hfill\\\\
(i) There exists a unique $4$-uple of scalars
$(\eta_1,\eta_2,\eta_3,\eta_4)$ such that the action by $g_1$,
$g_2$, $g_3$, $g_4$ on
$X=g_5v_3+\eta_1v_1+\eta_2v_2+\eta_3v_3+\eta_4v_4$ is a
multiplication by $r$.\\\\
(ii) For each integer $k\geq 6$, there exists a unique $4$-uple of
scalars $(\eta_1^k,\eta_2^k,\eta_3^k,\eta_4^k)$ such that the action
by $g_1$, $g_2$, $g_3$, $g_4$ on
$X_k=g_kv_3+\eta_1^kv_1+\eta_2^kv_2+\eta_3^kv_3+\eta_4^kv_4$ is a
multiplication by $r$.
\end{Claim}
\noin\textsc{Proof.} Immediate with $(\divideontimes)$ \\\\
Before we go further, we will need to have a better knowledge of
$v_4$, $v_2$ and $v_1$. We compute $v_4$ with the relation
$v_4=g_2v_3-r\,v_3$ and we get
\begin{multline}
v_4=\Bigg\lb\bigg(\unsur{l}-r\bigg)\bigg(\la_{12}+mr\,\wh{\la_{12}}-m\,\sum_{j=5}^n\la_{2j}r^{j-4}\bigg)-m\,(\la_{12}+\wh{\la_{12}})\Bigg\rb\,w_{12}\\-\wh{\la_{12}}(\wh{w_{24}}-w_{24})+r\,\wh{\la_{12}}(\wh{w_{14}}-w_{14})+\unsurr\,\wh{\la_{12}}(w_{23}-\wh{w_{23}})
+\wh{\la_{12}}(\wh{w_{13}}-w_{13})\\+\sum_{j=5}^n\la_{2j}(w_{1j}-\wh{w_{1j}})-\unsurr\,\sum_{j=5}^n\la_{2j}(w_{2j}-\wh{w_{2j}})
\end{multline}
Next, $v_2$ is given by formula $(56)$. In what follows, a term
carries a star if it is multiplied by a factor
$\unsur{r^2}+\unsur{r^4}$. So we have by $(56)$,
$$v_2^{^{\bigstar}}=-g_1v_3+\big(r-\unsur{r^3}-\unsurr\big)v_3-\unsur{r^2}v_4$$
We will study the coefficients of $w_{1,j}$ and $w_{2,j}$ for $j\geq
4$ in $v_2^{\star}$ and show that these coefficients are all zero.
Going back to the expression of $v_3$ in $(54)$, we see that the
coefficient of $w_{1j}$ in $-g_1v_3$ is
$(1+\frac{m}{r})\la_{2j}=\unsur{r^2}\la_{2j}$. This coefficient
cancels with the coefficient $-\unsur{r^2}\la_{2j}$ of $w_{1j}$ in
$-\unsur{r^2}v_4$. Similarly, the coefficient of $w_{14}$ in
$-g_1v_4$ is $-(m+r)\wh{\la_{12}}=-\unsurr\wh{\la_{12}}$, which
cancels to $\unsurr\wh{\la_{12}}$, the coefficient of $w_{14}$ in
$(r-\unsur{r^3}-\unsurr)\,v_3-\unsur{r^2}\,v_4$. When $j\geq 5$, the
coefficient of $w_{2j}$ in $v_2^{^{\star}}$ is
$(\unsurr-r)\,\la_{2j}+(r-\unsur{r^3}-\unsurr)\,\la_{2j}+\unsur{r^3}\,\la_{2j}$,
that is zero. And the coefficient of $w_{24}$ in $v_2^{^{\star}}$ is
$-mr\wh{\la_{12}}-(r-\unsur{r^3}-\unsurr)r\wh{\la_{12}}-\unsur{r^2}\wh{\la_{12}}$,
that is zero. Since $v_1=g_4v_2$, there are obviously no terms in
$w_{1,j}$ or $w_{2,j}$ for $j\geq 5$ in $v_1$ either. We are now in
a position to use the claim. First, by point $(ii)$ of Claim $2$, we
have for any $k\geq 6$
\begin{equation}
g_2X_k=r\,X_k
\end{equation}
Look at the coefficient of the term in $w_{2,k-1}$ in $(62)$. Using
the discussion above, it comes
$$r(\la_{2k}+\eta_3^k\,\la_{2,k-1}-\frac{\eta_4^k}{r}\,\la_{2,k-1})=-m\,\la_{2,k}-m\,\eta_3^k\,\la_{2,k-1}+\la_{2,k-1}\,\eta_4^k+\frac{m}{r}\,\la_{2,k-1}\,\eta_4^k,$$
from which we derive
\begin{equation}
\la_{2,k}=\la_{2,k-1}\,\big((r+\unsurr)\,\eta_4^k-\eta_3^k\big)
\end{equation}
Look now at the coefficient of the term in $w_{14}$ in the same
equation $(62)$. In order to get the coefficient on the left hand
side, we must in particular look at the coefficient in $w_{24}$ in
$g_4v_2$, so using the discussion above, we must look at the
coefficient in $w_{23}$ in $v_2$, as there is no term in $w_{24}$ in
$v_2$. Up to a division by a factor $\unsur{r^2}+\unsur{r^4}$, this
coefficient is
$$-\big(m(\la_{12}+\wh{\la_{12}})-\unsurr\,\wh{\la_{12}}+m\,\wh{\la_{12}}-m\,\sum_{j=5}^n\la_{2j}r^{j-5}\big)
-(r-\unsur{r^3}-\unsurr)(\la_{12}+\wh{\la_{12}})-\unsur{r^2}(\unsurr\,\wh{\la_{12}})$$
After replacing the sum as in $(58)$, all the terms simplify nicely
and yield $$\bigg(\unsurr+\unsur{r^3}\bigg)\la_{12}$$ Further the
coefficient of the term in $w_{14}$ in $g_4v_2$ is given by the
coefficient of the term in $w_{13}$ in $v_2$. Up to a division by a
factor $\unsur{r^2}+\unsur{r^4}$, this coefficient is
$$-\frac{m}{r}\la_{12}-\frac{m}{r}\wh{\la_{12}}+m\sum_{j=5}^nr^{j-6}\la_{2j}+\frac{\wh{\la_{12}}}{r}$$
By using again $(58)$ and simplifying, this coefficient is simply
$$\unsurr\bigg(\unsur{r}+\unsur{r^3}\bigg)\la_{12}$$
So, the coefficients of the term in $w_{14}$ in $\eta_1^k\,g_2\,v_1$
and in $r\,\eta_1^k\,v_1$ respectively cancel each other. We thus
obtain
$$-r^2\,\wh{\la_{12}}-r\,\wh{\la_{12}}\,\eta_3^k+\wh{\la_{12}}\,\eta_4^k=-r^2\,\wh{\la_{12}}\,\eta_4^k$$
Equivalently,
\begin{equation}
\big((1+r^2)\,\eta_4^k-r\,\eta_3^k\big)\,\wh{\la_{12}}=r^2\,\wh{\la_{12}}
\end{equation}
Assume for now that $\wh{\la_{12}}$ is nonzero. Then we get after
dividing also by $r$,
\begin{equation}
\forall k\geq 6,\;(r+\unsurr)\,\eta_4^k-\eta_3^k=r
\end{equation}
and so by $(63)$,
\begin{equation}
\forall k\geq 6,\,\la_{2k}=r\,\la_{2,k-1}
\end{equation}
We derive now from $(58)$
$$-\wh{\la_{12}}+\unsur{r^2}\,\la_{12}=m\,\la_{25}\,\sum_{j=5}^n
r^{j-6}r^{j-5}$$ And after evaluating the sum, it comes
\begin{equation}
-\wh{\la_{12}}+\unsur{r^2}\,\la_{12}=\bigg(\unsur{r^2}-r^{2n-10}\bigg)\,\la_{25}
\end{equation}
To get more relations between the coefficients, we use point $(i)$
of Claim $2$. First, we look at the coefficient of the term in
$w_{25}$ in $g_2X=r\,X$. We get after simplifications,
$$
r\,\la_{24}-mr\,\wh{\la_{12}}=\la_{25}\bigg(\frac{m}{r}-\frac{\eta_3}{r}+(1+\unsur{r^2})\eta_4\bigg)
$$
Recall from earlier that $\la_{24}=-\wh{\la_{24}}=-r\,\wh{\la_{12}}$
(see page $34$ of the present paper). Hence, we get
\begin{equation}
-\wh{\la_{12}}=\la_{25}\bigg(\unsur{r^2}-1-\frac{\eta_3}{r}+(1+\unsur{r^2})\eta_4\bigg)
\end{equation}
In particular, since we assumed that $\wh{\la_{12}}$ is non-zero, it
follows that $\la_{25}$ is also non-zero. Further, look at the
coefficient of the term in $w_{16}$ still in $g_2X=r\,X$ and derive
after using $(66)$ with $k=6$ and simplifying by $\la_{25}$
\begin{equation}
\eta_4=1
\end{equation}
Furthermore, look at the coefficient of the term in $w_{26}$ this
time in the same equation $g_2\,X=r\,X$. After simplifying and using
$(69)$, we get
\begin{equation}
\eta_3=\unsurr
\end{equation}
Plugging $(69)$ and $(70)$ in $(68)$ now yields
\begin{equation}
\la_{25}=-r^2\,\wh{\la_{12}}
\end{equation}
Next, by plugging $(71)$ into $(67)$, we obtain
\begin{equation}
\la_{12}=r^{2n-6}\,\wh{\la_{12}}
\end{equation}
Combining $(60)$ and $(72)$, we derive immediately
$$l=-\unsur{r^{2n-5}}$$
This is the value announced in Theorem $5$. It remains to show that
our assumption that $\wh{\la_{12}}\neq 0$ indeed holds. If
$\wh{\la_{12}}=0$, equation $(60)$ implies that $\la_{12}=0$. Then
many terms in $v_3$ and in $v_4$ vanish. Further, by $(68)$, we get
$\la_{25}=0$ or $\eta_3-(r+\unsurr)\,\eta_4=\unsurr-r$. Suppose the
second equality holds. Looking at the coefficient of the term in
$w_{2j}$, $j\geq 6$ in $g_2X=r\,X$ yields
$$\la_{2j}\bigg(1+\frac{\eta_3}{r}-(1+\unsur{r^2})\,\eta_4\bigg)=0,$$
and so $\unsur{r^2}\,\la_{2j}=0$. Then, $\la_{2j}=0$ for all $j\geq
6$. Next, look at the coefficient of the term in $w_{16}$ in
\begin{equation}g_2\,X_6=r\,X_6\end{equation} where we used the notations of
Claim $2$. Since $\la_{26}=0$, we simply get $\la_{25}=0$. But then
$v_3$ is zero, and this is a contradiction. So looking back up, we
must in fact have $\la_{25}=0$. We show this implies inductively
that all the $\la_{2j}$'s, $j\geq 6$ are zero. Let $k\geq 6$ and
suppose that $\la_{2,k-1}=0$. Let's show that $\la_{2,k}=0$ then. It
suffices to look at the coefficient of $w_{1,k}$ in
$$g_2\,X_k=r\,X_k$$
Get
$$-m\,\la_{2k}+\eta_3\,\la_{2k}-\unsurr\,\eta_4\,\la_{2k}=r\,\eta_4\la_{2k},$$
which after simplification rewrites
\begin{equation}
\big(m+(r+\unsurr)\,\eta_4-\eta_3\big)\,\la_{2k}=0
\end{equation}
As we have seen that $\eta_3-(r+\unsurr)\,\eta_4\neq m$, this forces
$\la_{2k}=0$, as announced. The fact that all the $\la_{2k},\,k\geq
5$ are zero is again a contradiction. So our initial assumption
$\wh{\la_{12}}=0$ is not possible. A consequence also is that
without loss of generality, $\wh{\la_{12}}$ can be set to $1$. Then
$\la_{12}$ is uniquely determined by $(60)$. And $\la_{25}$ is
uniquely determined by $(71)$. In turn, the $\la_{2k}$'s, $k\geq 6$
are uniquely determined by $(66)$. Thus, $v_3$ is uniquely
determined. And so is $v_4$. Then, $v_2$ is uniquely determined by
$(56)$ and $v_1$ is in turn completely determined by $v_1=g_4v_2$.
And so we have the following intermediate result.
\newtheorem{Result}{Result}
\begin{Result}
If in the Cohen-Wales space $V_n$, there exists an irreducible
invariant subspace that is isomorphic to $S^{(1),(n-1)}$, then
$l=-\unsur{r^{2n-5}}$ and there exists in $V_n$ a unique irreducible
$\mcalh(D_4)$-module that is isomorphic to $S^{(1),(3)}$.
\end{Result}

To finish the proof of Theorem $5$ stated in the introduction, it
remains to show that the Specht module $S^{(1),(1^{n-1})}$ cannot
occur in the Cohen-Wales space. Suppose it does, and let $\W$ be an
irreducible invariant subspace of $V_n$ that is isomorphic to
$S^{(1),(1^{n-1})}$. Then,
$$\W\da_{\mcalh(D_4)}\simeq (n-4)\,S^{(0),(1^4)}\oplus S^{(1),(1^3)}$$
We show that it is impossible to have vectors $y_1,y_2,y_3,y_4$ such
that the matrices of the left actions by the $g_k$'s with $k\in\lb
1,2,3,4\rb$ on these vectors is given by the matrices $H_i$'s above,
where $r$ has been replaced by $-\unsurr$. Let's call these
conjugate matrices the $K_i$'s. First, the set of relations
$$\left\lb\begin{array}{ccc}
g_2y_1&=&-\unsurr\,y_1\\
g_3y_1&=&-\unsurr\,y_1
\end{array}\right.$$
forces without loss of generality
\begin{equation}y_1=w_{13}-\unsurr\,w_{23}-r\,w_{12}\end{equation}
From there, it is very easy to conclude. Indeed, $y_2$ is given by
the first column of $K_4$, then $y_3$ is provided by the second
column of $K_3$ and finally $y_4$ is given by third column of $K_2$.
Since an action by the $g_i$'s with $i=2,3,4$ on "non-hat terms"
will never create a "hat term" by defining equations
$(8),(9),(11),(14)$, $(15),(16)$ of the representation, we see with
$(75)$ that the $y_i$'s do not contain any hats. However, an action
by $g_1$ on $y_1$ creates a term in $\wh{w_{23}}$ with coefficient
$1$. So $g_1y_1$ cannot be a linear combination of $y_1,y_2,y_3$ and
$y_4$. This is a contradiction. Hence we are done with the complete
proof of Theorem $5$, points $(i)$ and $(ii)$.



\subsection{Proof of the Main Theorem}
\subsubsection{Proof of the necessary condition}
In this part, we assume that $\ih(n)$ and $\mcalh(D_n)$ are
semisimple and we show that if $\n^{(n)}$ is reducible, then
$l\in\lb\unsur{r^{4n-7}},\unsur{r^{2n-7}},-\unsur{r^{2n-5}},r^3,-r^3,\unsurr\rb$.
We solve the small cases $n\in\lb 4,5,6,7,8,9\rb$ by computer means
when $n\in\lb 4,5,6,7\rb$ and by hand when $n\in\lb 8,9\rb$. We
explain them later. For now suppose the necessary condition above
holds in these cases and fix $n\geq 10$. We proceed by induction. We
assume that the necessary condition above holds at ranks $n-2$ and
$n-1$ and show it then holds at rank $n$. Suppose $\n^{(n)}$ is
reducible and let $\W$ be an irreducible proper invariant subspace
of $V_n$. Suppose
$l\not\in\lb\unsur{r^{4n-7}},\unsur{r^{2n-7}},-\unsur{r^{2n-5}}\rb$.
Then, by Theorems $3,4,5$ and Theorem $8$, point $(i)$, we must have
$$\text{dim}(\W)\geq\cil$$
When $n\geq 10$, we claim that the dimension of $\W$ is large enough
so that the intersection spaces $\W\cap V_{n-2}$ and $\W\cap
V_{n-1}$ are nonzero. Indeed, if $\W\oplus V_{n-2}$, then it comes
$\text{dim}(\W)\leq \text{dim}(V_n)-\text{dim}(V_{n-2})=4n-6$. And
if $\W\oplus V_{n-1}$, then it comes $\text{dim}(\W)\leq
\text{dim}(V_n)-\text{dim}(V_{n-1})=2n-2$. But
$$\cil>4n-6\Longleftrightarrow n\geq 10\;\;\text{and}\;\;\cil>2n-2\Longleftrightarrow n\geq
7$$ So when $n\geq 10$, both intersections are nonzero. Moreover,
both intersections are proper, because if $\W$ contains $V_{n-1}$ or
$V_{n-2}$, then it is quite easy to see on the representation that
$\W$ would be the whole space $V_n$. From there, we get that both
$\n^{(n-1)}$ and $\n^{(n-2)}$ are reducible, so applying the
induction hypothesis yields
$$\large{\left\lb\begin{array}{l}
l\in\lb\unsur{r^{4n-11}},\unsur{r^{2n-9}},-\unsur{r^{2n-7}},r^3,-r^3,\unsurr\rb\\
\&\\
l\in\lb\unsur{r^{4n-15}},\unsur{r^{2n-11}},-\unsur{r^{2n-9}},r^3,-r^3,\unsurr\rb
\end{array}\right.}$$

\noindent With our restrictions on $r$, we see that $l$ must then
belong to the set of values $\lb r^3,-r^3,\unsurr\rb$.
This finishes the proof in the general case.\\

Let's now deal with the case $n=9$, still assuming the result holds
for the smaller values of $n$. Suppose $\n^{(9)}$ is reducible and
let $\W$ be an irreducible invariant subspace of $V_9$. Using part
$\S 3.2$, the degrees less than $72$ of the irreps of $\mcalh(D_9)$
are
$$1,8,9,27,28,36,42,48,56,63,70$$ If
$l\not\in\lb\unsur{r^{29}},\unsur{r^{11}},-\unsur{r^{13}}\rb$, then
by Theorems $3,4,5$ we must have $\text{dim}(\W)\geq 27$. First, if
$\text{dim}(\W)\geq 36>2(18-3)=30$, then $\W\cap V_7\neq\lb 0\rb$
and the general technique applies. Hence suppose
$\text{dim}(\W)\in\lb 27,28\rb$. Define $\W_8=\W\cap V_8$. We have
$$\text{dim}(\W_8)\geq\text{dim}(\W)+\text{dim}(V_8)-\text{dim}(V_9)$$
So, if $\text{dim}(\W)=27$, we get $\text{dim}(\W_8)\geq 11$ and if
$\text{dim}(\W)=28$, we get $\text{dim}(\W_8)\geq 12$. In any case,
if $\W_8$ is irreducible, we must have $\text{dim}(\W_8)\geq 14$ by
Theorem $8$, point $(iii)$. Assume first that $\W_8$ is reducible.
If $\text{dim}(\W_8)>14=2(8-1)$, then $\n^{(7)}$ is reducible. Then
both $\n^{(8)}$ and $\n^{(7)}$ are reducible and we conclude like in
the general case. So, suppose $\text{dim}(\W_8)=14$. Then notice
$\W_8\oplus V_7=V_8$ by an inspection on the dimensions. In
particular, there exists elements $z_1\in\W_8$ and $z_2\in V_7$ such
that
$$w_{78}=z_1+z_2$$
It then comes
$$w_{89}=e_9\,.\,w_{78}=e_9\,.\,z_1\;\;,$$
as the tangle resulting from an action by $e_9$ on any basis vector
of $V_7$ contains two horizontal strands at the bottom: one joining
nodes $6$ and $7$, the other one joining nodes $8$ and $9$. By
construction of the representation, such an element is zero. Now,
$z_1$ belongs to $\W$. So, $e_9\,.\,z_1$ belongs to $\W$. Then,
$w_{89}$ belongs to $\W$. This implies that $\W$ is the whole space
$V_9$, a contradiction. \\
Suppose now $\W_8$ is reducible. By the semisimplicity assumption
for $\mcalh(D_8)$, the fact that $\text{dim}(\W_8)\geq 11$ and the
uniqueness part in Theorem $3$, we must again have
$\text{dim}(\W_8)\geq 14$. So this case reduces to the previous
case. This ends the case $n=9$.\\

Let's now deal with the cases $n\in\lb 4,5,6\rb$. Suppose $\n^{(n)}$
is reducible and let $\W$ be such a proper invariant subspace of
$V_n$. By Proposition $1$ of $\S 3.1$, we know that
\begin{equation}e_k.\W=0\qquad\forall 1\leq k\leq n\end{equation}
\newtheorem{Definition}{Definition}
\begin{Definition}
We define algebra elements that are some conjugates of the $e_k$'s.
\begin{eqnarray*}
C_{i,i+1}&=&e_{i+1}\qquad\qquad\qquad\;\;\;\forall 1\leq i\leq n-1\\
C_{i,j}&=&g_{j,i+2}\,e_{i+1}\,g_{i+2,j}^{\star}\qquad\forall 1\leq i<j\leq n\;\text{with}\;j\geq i+2\\
\wh{C_{12}}&=&e_1\\
\wh{C_{1,j}}&=&g_{j,3}\,e_1\,g_{3,j}^{\star}\qquad\qquad\;\;\;\forall 3\leq j\leq n\\
\wh{C_{i,j}}&=&g_{i,2}\,g_{j,3}\,e_1\,g_{3,j}^{\star}\,g_{2,i}^{\star}\qquad\forall
2\leq i<j\leq n
\end{eqnarray*}
By $g_{s,t}$ (resp $g_{s,t}^{\star}$), we understand the product of
the $g_k$'s (resp the $g_k^{-1}$s), where $k$ lies on the integer
path from $s$ up or down to $t$.
\end{Definition}
\begin{center}
\epsfig{file=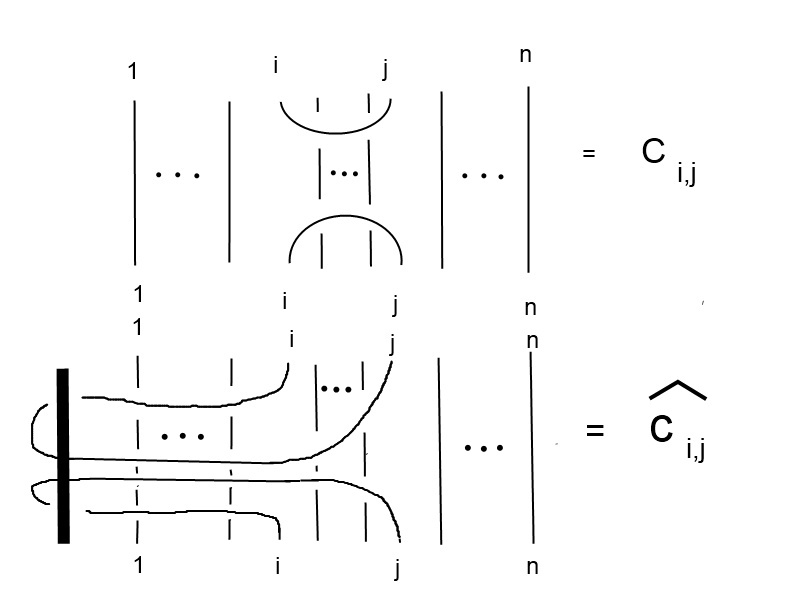, height=9cm}
\end{center}
\begin{Definition}
Define $$S(n)=\sum_{1\leq i<j\leq n}C_{ij}+\sum_{1\leq i<j\leq
n}\wh{C_{ij}}$$
\end{Definition}
\begin{Definition}
Define $$K(n)=\bigg(\cap_{1\leq i<j\leq
n}\;\text{Ker}\;\n^{(n)}(C_{ij})\bigg)\cap\bigg(\cap_{1\leq i<j\leq
n}\;\text{Ker}\;\n^{(n)}(\wh{C_{ij}})\bigg)$$ and let $k(n)$ be the
dimension of $K(n)$ as a vector space over $\Q(l,r)$.
\end{Definition}
After this series of definitions, we are back to the proof of the
necessary condition. Equalities $(76)$ and the fact that $\W$ is
invariant imply that the left action by $S(n)$ on $\W$ is trivial.
Then, since $\W\neq\lb 0\rb$, the determinant of this action must be
zero. Using the tangles, we computed the matrix of the left action
by $S(n)$ in the basis formed by the $w_{ij}$'s and the
$\wh{w_{ij}}$'s. Note each row of this matrix corresponds to the
action of one of the $C_{ij}$'s or $\wh{C_{ij}}$'s. In particular,
ordering the vectors of the basis in such a way that the $2(n-1)$
last vectors of this basis have an extremity ending in $n$ allows us
to only have to compute $4(n-1)^3$ entries and not $n^2(n-1)^2$
after rank $4$. For $n=4$, the matrix is of size $12$, when $n=5$ of
size $20$ and when $n=6$ of size $30$. We used Mathematica to solve
the determinant of this matrix equals zero. We obtained the values
of Theorem $2$ point $(i)$ in each of the cases $n=4,5,6$. When
$n=7$, the number of entries is unreasonably big to do it by hand.
Thus, we wrote a program in Mathematica that computes the sum
matrix. All the matrices in that program are defined by blocks and
inductively.
When running the program, we obtained
$$\text{det}\,S(7)=\frac{(-1+lr)^{21}(l-r^3)^{14}(l+r^3)^{35}(-1+lr^7)^6(1+lr^9)^7(-1+lr^{21})}{l^{42}\,r^{105}\,(r^2-1)^{42}}$$
So, if $\n^{(7)}$ is reducible, then
$$l\in\bigg\lb\unsur{r^{21}},\unsur{r^7},-\unsur{r^9},r^3,-r^3,\unsurr\bigg\rb\qquad\qquad(\clubsuit)$$
We get the values of Theorem $2$ point $(i)$ for $n=7$. This terminates the case $n=7$.\\

Let's finish all the cases by doing $n=8$. Suppose $\n^{(8)}$ is
reducible and let $\W$ be an irreducible invariant subspace of
$V_8$. By the discussion of $\S3.2$ of this paper and Appendix B of
\cite{THE}, the degrees of the irreps of $\mcalh(D_8)$ that are less
than $56$ are
$$1,7,8,14,20,21,28,35,42,48$$ Suppose $l\not\in\lb\unsur{r^{25}},\unsur{r^9},-\unsur{r^{11}}\rb$.
Let's show that $l\in\lb\unsurr,r^3,-r^3\rb$ then. By Theorems
$3,4,5$, we have $\text{dim}(\W)\geq 14$. Further, we have
$2(2\times 8-3)=26$, so if $\text{dim}(\W)\geq 28$, then the general
technique applies. \\ - Hence suppose $\text{dim}(\W)=21$ and so
$\W$ is isomorphic to $S^{(0),(6,1,1)}$ or its conjugate
$S^{(0),(3,1^5)}$. Then,
$$\W\da_{\mcalh(D_7)}\simeq S^{(0),(5,1,1)}\oplus S^{(0),(6,1)}\;\;\text{or}\;\;
\W\da_{\mcalh(D_7)}\simeq S^{(0),(2,1^5)}\oplus S^{(0),(3,1^4)}$$
Now look at $\W\cap V_7$. We have $$\text{dim}(\W\cap V_7)\geq
21+42-56=7$$ So $\W\cap V_7$ cannot be isomorphic to $S^{(0),(6,1)}$
or to $S^{(0),(2,1^5)}$. Further, $\W\cap V_7$ is not $\W$ either,
as otherwise $\W$ would be the whole space $V_8$. Then $\W\cap V_7$
must be isomorphic to $S^{(0),(5,1,1)}$ or to $S^{(0),(3,1^4)}$. In
any case, the dimension of $\W\cap V_7$ is $15$. Moreover, since
$\W\cap V_7$ is a proper invariant subspace of $V_7$, it is
annihilated by all the $\ovl{C_{ij}}$'s with $1\leq i<j\leq 7$.
Hence we have $\W\cap V_7\subseteq K(7)$, and so $k(7)\geq 15$.
\begin{Lemma}
$$k(n)=n^2-n-\text{rank}(S(n))\qquad\forall n\geq 4$$
\end{Lemma}
\noindent\textsc{Proof of the lemma.} Obvious by the remark above
that each row of the matrix of the left action by $S(n)$ on the
basis vectors $\ovl{w_{ij}}$'s, $1\leq i<j\leq n$ corresponds to the
action of one of the $\ovl{C_{ij}}$'s, as the kernel of the sum
matrix is then $K(n)$. \\
We computed the rank of $S(7)$ with Mathematica for the different
values of $l$ and $r$ present in $(\clubsuit)$. Here is what we got.
$$\begin{array}{ccc}
\text{When}&l=\unsur{r^{21}},&k(7)=1\\
&&\\
\text{When}&l=\unsur{r^7},&k(7)=6\\
&&\\\text{When}&l=-\unsur{r^9},&k(7)=7\\
&&\\\text{When}&l=r^3,&k(7)=14\\
&&\\\text{When}&l=\unsurr,&k(7)=21\\
&&\\\text{When}&l=-r^3,&k(7)=35\\
\end{array}$$
So, if $k(7)\geq 15$, this forces $l\in\lb\unsurr,-r^3\rb$.\\\\
- Suppose now $\text{dim}(\W)=20$. By the discussion of $\S\,3.2$
and the table of Appendix $B$ of \cite{THE}, the only classes of
irreducible $\mcalh(D_8)$-modules of dimension $20$ are
$S^{(0),(6,2)}$ and $S^{(0),(2^2,1^4)}$. Then, by Lemma $3$, we get
$\W\simeq S^{(0),(6,2)}$ and $l=r^3$.\\
- Suppose finally $\text{dim}(\W)=14$. Then $\W\simeq S^{(0),(4,4)}$
or $\W\simeq S^{(0),(2,2,2,2)}$.
\begin{itemize}
\item If $\W\simeq S^{(0),(4,4)}$, then $\W\da_{\mcalh(D_7)}\simeq S^{(0),(4,3)}$
and so $\W\da_{\mcalh(D_6)}\simeq S^{(0),(3,3)}\oplus
S^{(0),(4,2)}$. Then the Specht module $S^{(0),(3,2)}$ is a
constituent of $\W\da_{\mcalh(D_5)}$. Then, the proof of point $(i)$
of Lemma $3$ shows that $l$ must be equal to $r^3$.
\item If $\W\simeq S^{(0),(2^4)}$, then $\W\da_{\mcalh(D_7)}\simeq
S^{(0),(2^3,1)}$ and so $\W\da_{\mcalh(D_6)}\simeq
S^{(0),(2^3)}\oplus S^{(0),(2,2,1,1)}$. Then the Specht module
$S^{(0),(2,2,1)}$ is a constituent of $\W\da_{\mcalh(D_5)}$. The
proof of point $(ii)$ of Lemma $3$ shows that this cannot happen.
\end{itemize}
So, we are done with the proof of the necessary condition.

\subsubsection{Proof of the sufficient condition}

The reducibility of $\n^{(n)}$ is already known when
$l=\unsur{r^{4n-7}}$ by Theorem $3$ and when $l=\unsur{r^{2n-7}}$ by
Theorem $4$. Thus, it remains to show that the representation is
reducible when $l\in\lb -\unsur{r^{2n-5}},r^3,-r^3,\unsurr\rb$.
First when $l=\unsurr$, we read on the representation that the
vectors $t_{ij}$'s of Theorem $7$ span an
$\frac{n(n-1)}{2}$-dimensional invariant subspace of $V_n$ and so
the representation is reducible in this case as well. When $l\in\lb
r^3,-r^3,-\unsur{r^{2n-5}}\rb$, we need a lemma to show the
reducibility of the representation.
\begin{Lemma} Let $n$ be an integer with $n\geq 4$.
The vector space $K(n)$ is a $\cgwn$-submodule of $V_n$.
\end{Lemma}
\textsc{Proof of the Lemma.} We want to show that
$$x\in\cap_{1\leq i<j\leq n}\,Ker\,\ovl{C_{ij}}\Rightarrow g_k\,x\in\cap_{1\leq i<j\leq
n}\,Ker\,\ovl{C_{ij}}\qquad\forall 1\leq k\leq n$$ Let $x\in K(n)$.
We proceed in four steps. Step $1$ and Step $2$ are almost
identical, Step $3$
uses Step $2$ and Step $4$ uses Step $1$ and Step $3$.\\
- \textbf{Step $1$}. We show that $g_k\,x\in\cap_{1\leq i<j\leq
n}\,Ker\,C_{ij}$ for every integer $k$ with $2\leq k\leq n$. Fix
integers $i$ and $j$ with $1\leq i<j\leq n$ and $j\geq i+2$. Let's
show that $C_{ij}g_kx=0$.
\\First, if $k\in\lb i,i+1,j,j+1\rb$, acting to the right of $C_{ij}$
with $g_k$ shifts one of the extremities of the bottom horizontal
strand (use the Kauffman skein relation when necessary) and the
result follows. When $j=i+1$, use the delooping relation for the
"middle case". \\
Next, if $i+2\leq k\leq j-1$, use Reidemeister's move
\Roman{chiromains} twice to notice that $C_{ij}g_k=g_kC_{ij}$, and
so $C_{ij}g_kx=0$.\\ Finally, when $k\leq i-1$ or $k\geq j+2$, the
elements $C_{ij}$ and $g_k$ also commute, which gives the result
in this case as well.\\
- \textbf{Step $2$}. We show that $g_kx\in\cap_{1\leq i<j\leq
n}\,Ker\,\wh{C_{ij}}$ for all $k$ with $2\leq k\leq n$. This case is
identical as Step $1$, except Reidemeister's move \Roman{chiromains}
must also be used twice when $k\leq i-1$ and when $j=i+1$, we use
$$\wh{C_{i,i+1}}\,g_{i+1}\,x=\delta^{-1}\,\Xi^{+}\,\wh{C_{i,i+1}}\,x,$$
after moving the crossing near the pole thanks to Reidemeister's
move \Roman{chiromains}. \\
- \textbf{Step $3$}. We show that $g_1x\in\cap_{1\leq i<j\leq
n}\,Ker\, C_{ij}$. Because $g_1$ commutes to $C_{ij}$ when $i\geq
3$, we only need to show that $C_{2,j}g_1x=0$ and $C_{1,j}g_1x=0$.
But,
$$C_{2,j}g_1=g_1g_1^{-1}C_{2j}g_1=g_1\wh{C_{1j}},$$ where the second
equality follows from an application of the double twist relation,
as in Figure $3$. The second one is more difficult and requires
careful manipulations on the tangles. Proceed as follows. Multiply
the top of $C_{1j}$ by $g_1g_1^{-1}$. Use the commuting relation at
the top and at the bottom, as on Figure $4$. Then use the double
twist relation and get all together $g_1$ times a tangle that is
almost $\wh{C_{2j}}$, except the top horizontal strand and the
bottom horizontal strand both under-cross the vertical strand that
they first intersect when sliding along the strands from the left
hand side extremities. Now, it suffices to multiply at the top by
$g_2$ and at the bottom by $g_2^{-1}$ to get $\wh{C_{1,j}}$ instead.
So, we have
$$g_1g_1^{-1}C_{1,j}g_1x=g_1(g_2^{-1}\wh{C_{1j}}g_2)x$$
But by Step $2$, we have $\wh{C_{1j}}\,g_2\,x=0$, so we are done.
\vspace{-0.6cm}\begin{center}\epsfig{file=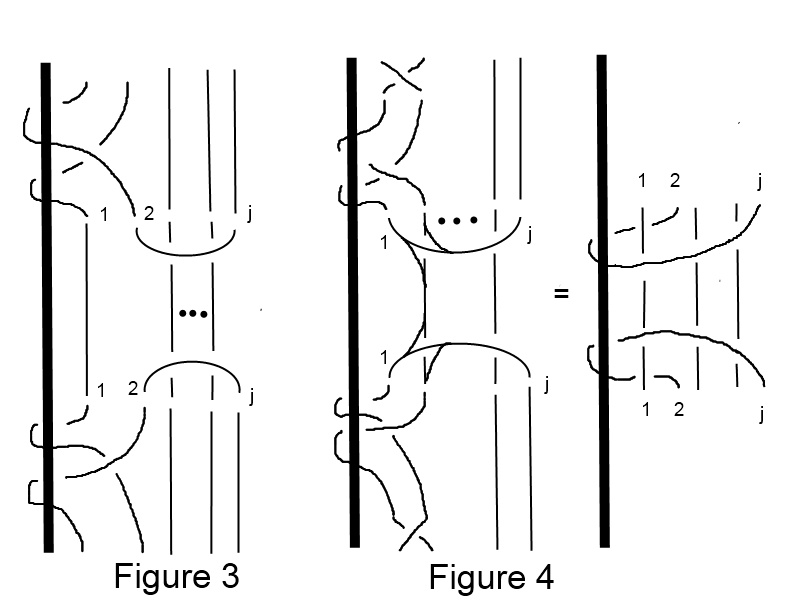,
height=10cm}\end{center} - Step $4$. We show that
$g_1x\in\cap_{1\leq i<j\leq n}\,Ker\,\wh{C_{ij}}$. First, we deal
with $\wh{C_{1j}}$ and $\wh{C_{2j}}$. With the Kauffman skein
relation, it suffices to show that $\wh{C_{1j}}g_1^{-1}x=0$ and
$\wh{C_{2j}}g_1^{-1}x=0$. We have
$$\begin{array}{ccccccc}
\wh{C_{1j}}g_1^{-1}x&=&g_1^{-1}g_1\wh{C_{1j}}g_1^{-1}x&\overset{1}{=}&g_1^{-1}C_{2j}x&=&0\\
\wh{C_{2j}}g_1^{-1}x&\overset{2}{=}&g_1^{-1}g_2C_{2j}g_2^{-1}x&\overset{3}{=}&0&&
\end{array}$$
Equality $1$ is obtained by using in the respective order,
simultaneously at the top and at the bottom the double twist
relation followed by a Reidemeister's move
\setcounter{chiromains}{2} \Roman{chiromains}, then the double twist
relation again.\\
For equality $2$, multiply to the left by $g_1^{-1}g_1$, use the
commuting relation at the top and at the bottom, followed by the
double twist relation and a Reidemeister's move \Roman{chiromains}.
Get $g_1^{-1}$ times the tangle of Figure $5$ below. The latter
tangle is $g_2C_{1j}g_2^{-1}$, hence equality $2$.\\ Equality $3$ is
obtained by using the Kauffman skein relation and Step $1$.\\
It remains to show that $\wh{C_{ij}}g_1x=0$ when $3\leq i<j\leq n$,
or which is equivalent $\wh{C_{ij}}g_1^{-1}x=0$. We chose to do it
algebraically. When $i\geq 4$, we have
$$\wh{C_{ij}}=g_{i,4}\,\wh{C_{3j}}\,g_{4,i}^{\star}$$
We compute $\wh{C_{3,j}}\,g_1^{-1}$. We have, where the parenthesis
point out where the next transformations take place.
\begin{eqnarray}\wh{C_{3,j}}\,g_1^{-1}&=&g_3\,g_2\,g_{j,3}\,e_1\,g_{3,j}^{\star}\,g_2^{-1}(g_3^{-1}g_1^{-1}g_3^{-1})g_3\\
&=&g_3\,g_2\,g_{j,3}\,e_1\,g_{3,j}^{\star}\,(g_2^{-1}g_1^{-1})g_3^{-1}g_1^{-1}g_3\\
&=&g_3\,g_2\,(g_{j,4}\,e_3\,g_{4,j}^{\star})\,g_2^{-1}g_3^{-1}g_1^{-1}\\
&=&g_3\,g_2\,C_{2,j}\,(g_2^{-1})g_3^{-1}g_1^{-1}\\
&=&\!\!g_3g_2C_{2,j}g_2g_3^{-1}g_1^{-1}\!\!\!+\!mg_3g_2C_{2,j}g_3^{-1}g_1^{-1}\!\!-\!mg_3g_2C_{2,j}e_2g_3^{-1}g_1^{-1}
\end{eqnarray}
Equality $(78)$ is obtained by using the braid relation. To get
$(79)$, commute $g_1^{-1}$ to the right of $g_3^{-1}$, add a factor
$g_1g_1^{-1}$ in between $e_1$ and $g_3^{-1}$, use the braid
relation with nodes $1$ and $3$ and use the first delooping relation
$(DL)$. Now derive from the first equality in $(6)$ of Proposition
$2.3$ of \cite{CGW} that
\begin{equation}e_1g_3^{-1}g_1^{-1}=g_3^{-1}g_1^{-1}e_3\end{equation}
Further cancel the product $g_3g_3^{-1}$ to the left of $e_3$ and
cancel the same product to the right of $e_3$ after replacing $e_3$
with $l\,e_3g_3$. Commute $g_1^{-1}$ to the right hand side and use
the braid relation with nodes $1$ and $3$. Cancel the product
$g_3^{-1}g_3$ of the extreme left. Get $(79)$.\\
Equality $(81)$ is obtained by applying the Kauffman skein relation.
We study the three terms of this sum separately. Let's call them
$a$, $b$ and $c$. We have
$$a=g_3g_2^2C_{1j}g_3^{-1}g_1^{-1}=g_3g_2^2g_3^{-1}C_{1j}g_1^{-1}$$
Now the fact that $a$ annihilates $x$ follows from Step $3$.
Further, we have
$$b=mg_3g_2g_3^{-1}C_{3,j}g_1^{-1}=mg_3g_2g_3^{-1}g_1^{-1}C_{3,j}$$ and so
$b$ also annihilates $x$. Finally, we have
$$c=-mg_3g_2C_{2j}g_3^{-1}g_3e_2g_3g_1^{-1}=-mg_3g_2C_{2j}g_3^{-1}C_{13}g_1^{-1}$$
and again the fact that $c$ annihilates $x$ follows from Step $3$.\\
This settles Lemma $5$. The goal next is to show that this submodule
of $V_n$ is non-zero when $l\in\lb -\unsur{r^{2n-5}},r^3,-r^3\rb$.
The results are summarized in the following Theorem.
\begin{theo} (Reducibility of the representation $\n^{(n)}$ when $l\in\lb
r^3,-r^3,-\unsur{r^{2n-5}}\rb$).\\\\
(i) When $l=r^3$, the vector
$\mx=(w_{24}+r^2\,\wh{w_{24}})-r(w_{14}+r^2\wh{w_{14}})-r\,(w_{23}+r^2\,\wh{w_{23}})+r^2\,(w_{13}+r^2\wh{w_{13}})$
belongs to $K(n)$ for all $n\geq 4$.\\
(ii) When $l=-r^3$, the vector
$\my=w_{34}-\unsurr\,w_{35}+\unsur{r^2}w_{45}$ belongs to $K(n)$ for
all $n\geq 5$. The vector $\mz=r^3\wh{w_{24}}-r^2\wh{w_{34}}+w_{23}$
belongs to K(4).\\
(iii) When $l=-\unsur{r^{2n-5}}$, the vector
\begin{multline*}\mathcal{J}_n=
(\wh{w_{12}}+r^{2n-6}\,w_{12})-\unsurr\,(\wh{w_{13}}+r^{2n-6}\,w_{13})-(1+r^{2n-6})\,w_{23}\\+r\,\sum_{j=4}^nr^{j-5}\big\lb(w_{3,j}-\wh{w_{3,j}})-r\,(w_{2,j}-\wh{w_{2,j}})\big\rb
\end{multline*}
belongs to $K(n)$ for all $n\geq 4$.
\end{theo}
\noin Note $\mx$ is up to a sign the vector $v_4$ of Lemma $3$ in
$\S\,3.4$. Also, $\my$ is the vector $v_1$ of expression $(29)$ of
$\S\,3.4$. As for $\mz$ it was found with Mathematica. Finally,
$\mathcal{J}_n$ is the vector $v_3$ of expression $(54)$ of
$\S\,3.5$, where $\wh{\la_{12}}$ has been set to the value $1$ and
where the other coefficients are given by $(66)$, $(71)$ and $(72)$.
\vspace{-0.5cm}
\begin{center}
\epsfig{file=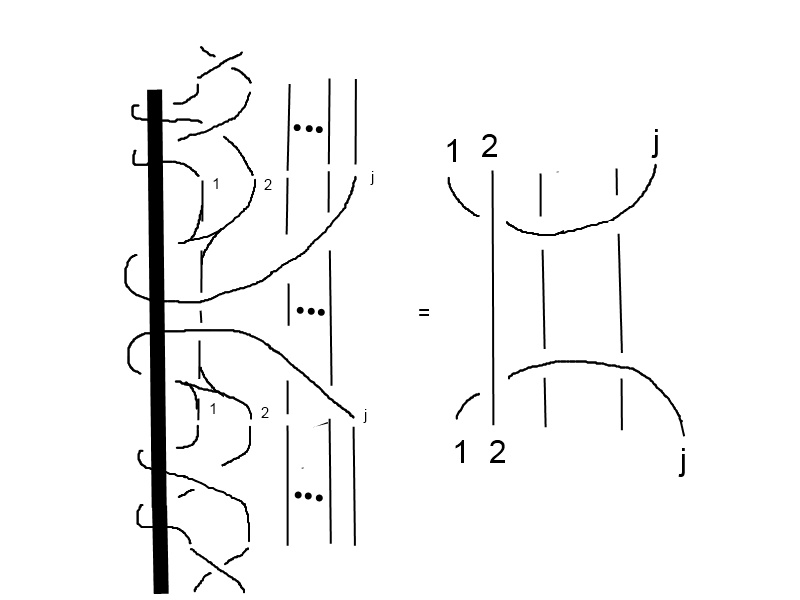, height=10cm}
\end{center}
\begin{center}
\textbf{\Large{Figure $5$}}
\end{center}
To prove Theorem $10$, we will make an extensive use of the
following proposition that provides the action by the $\ovl{C_{ij}}$
on the basis vectors of the Cohen-Wales space.
\begin{Proposition} The following equalities hold.
\begin{eqnarray*}
C_{ij}.\,w_{i-s,i}&=&\unsur{r^{(j-i)+(s-2)}}\,\,w_{ij}\;\;\qquad LONHNH_{j-i,s}\\
C_{ij}.\,\wh{w_{i-s,i}}&=&\unsur{r^{(j-i)+(s-2)}}\,\,w_{ij}\;\;\qquad LONHH_{j-i,s}\\
\wh{C_{ij}}.\,w_{i-s,i}&=&\frac{l}{r^{(j-i)+(s-3)}}\,\,\wh{w_{ij}}\;\;\qquad LOHNH_{j-i,s,l}\\
\wh{C_{ij}}.\,\wh{w_{i-s,i}}&=&\unsur{l\,r^{(j-i)+(s-1)}}\,\,\wh{w_{ij}}\;\;\qquad
LOHH_{j-i,s,l}\end{eqnarray*}
\begin{eqnarray*}
C_{ij}.\,w_{i,j-s}&=&\unsur{l\,r^{s-1}}\,\,w_{ij}\;\;\qquad LINHNH_{s,l}\\
C_{ij}.\,\wh{w_{i,j-s}}&=&\unsur{r^{s-2}}\,\,w_{ij}\;\;\qquad LINHH_{s}\\
\wh{C_{ij}}.\,\wh{w_{i,j-s}}&=&\unsur{l\,r^{s-1}}\,\,\wh{w_{ij}}\;\;\qquad LIHH_{s,l}\\
\wh{C_{ij}}.\,w_{i,j-s}&=&\unsur{r^{s-2}}\,\,\wh{w_{ij}}\;\;\qquad
LIHNH_{s}
\end{eqnarray*}
\begin{eqnarray*}
C_{ij}.\,\wh{w_{i+t,j-s}}&=&\frac{m\,r^{t-s-2}}{l}\,(1-l\,r)(1+r^2)\,w_{ij}\qquad\;\;\; INHH_{l,t-s}\\
\wh{C_{ij}}.\,w_{i+t,j-s}&=&0\qquad\qquad\qquad\qquad\qquad\qquad\qquad\qquad IHNH\\
\wh{C_{ij}}.\,\wh{w_{i+t,j-s}}&=&0\qquad\qquad\qquad\qquad\qquad\qquad\qquad\qquad IHH\\
\wh{C_{ij}}.\,\wh{w_{i-s,i-t}}&=&\frac{m}{l\,r^{j-i+s+t-2}}(1-lr)(1+r^2)\,\,\wh{w_{ij}}\qquad
ELOHH_{l,j-i,t+s}
\end{eqnarray*}
\begin{eqnarray*}
C_{ij}.\,w_{i-s,j-t}&=&\frac{m}{r^{s+t-2}}\bigg(\unsur{l}-\unsurr\bigg)\,\,w_{ij}\qquad\qquad LCNHNH_{l,s+t}\\
C_{ij}.\,\wh{w_{i-s,j-t}}&=&\frac{m}{r^{s+t-2}}\bigg(\unsur{l}-\unsurr\bigg)\,\,w_{ij}\qquad\qquad LCNHH_{l,s+t}\\
\wh{C_{ij}}.\,w_{i-s,j-t}&=&\frac{m}{r^{t+s-2}}(r-l)\,\,\wh{w_{ij}}\qquad\qquad LCHNH_{l,t+s}\\
\wh{C_{ij}}.\,\wh{w_{i-s,j-t}}&=&\frac{m}{r^{t+s-2}}\bigg(\unsur{l}-\unsurr\bigg)\,\wh{w_{ij}}\qquad\qquad
LCHH_{l,t+s}
\end{eqnarray*}
\end{Proposition}
In this proposition, the capital letters stand for the following words. \\
L: left; I: inside; H: hat; NH: non hat; E: extreme; O: outside; C:
crossed.\\\\
All these equalities were obained by using the tangles. For now
assume that they hold and let's prove the Theorem. Let's deal with
$(i)$. The program of the Appendix provides what we called the sum
matrix. Running it for $n=4$ and for $n=5$, we can check that these
matrices both annihilate $\mx$ and so $\mx$ lies in the intersection
$K(4)\cap K(5)$. For larger $n$, we proceed by induction. Let $n\geq
6$ and suppose that $\mx\in K(n-1)$. We must study the action by
$\ovl{C_{k,n}}$ on the vectors $\ovl{w_{13}}$, $\ovl{w_{23}}$,
$\ovl{w_{14}}$ and $\ovl{w_{24}}$. First, when $5\leq k\leq n-1$,
the action by $C_{k,n}$ on these vectors is zero. We next deal with
the actions by $\ovl{C_{1,n}}$, $\ovl{C_{2,n}}$, $\ovl{C_{3,n}}$ and
$\ovl{C_{4,n}}$. We see with the second set of formulas above that
$\ovl{C_{1,n}}(-r(w_{14}+r^2\,\wh{w_{14}})+r^2(w_{13}+r^2\,\wh{w_{13}}))=0$
and this independently from the values of $l$ and $r$. Moreover, the
first equality of the third set above implies that
$C_{1,n}.\,(r^2\,\wh{w_{24}}-r^3\,\wh{w_{23}})=0$. And so,
$C_{1,n}.\,\mx=0$. Also, by IHNH and IHH, we have
$\wh{C_{1,n}}.(w_{24}+r^2\wh{w_{24}}-r(w_{23}+r^2\wh{w_{23}}))=0$.
Hence, $\wh{C_{1,n}}.\,\mx=0$. For the action by $\ovl{C_{2,n}}$,
notice that \begin{eqnarray*}
\ovl{C_{2,n}}&=&g_2\,\ovl{C_{1,n}}\,g_2^{-1}\\
\text{and}\;\;g_2.\,\mx&=&-\unsurr\mx\;\;\text{when $l=r^3$}
\end{eqnarray*}
Hence $\ovl{C_{2,n}}.\,\mx=0$ by the previous case. By the first set
of relations above, the action by $\ovl{C_{3,n}}$ on the linear
combination of the vectors ending in node $3$ in $\mx$ is zero. And
by the last set of relations above, the action by $\ovl{C_{3,n}}$ on
the rest of $\mx$ is also zero. Hence $\ovl{C_{3,n}}.\,\mx=0$. This
also implies that $\ovl{C_{4,n}}.\,\mx=0$ after noticing that
\begin{eqnarray*}
\ovl{C_{4,n}}&=&g_4\,\ovl{C_{3,n}}\,g_4^{-1}\\
\text{and}\;\;g_4.\,\mx&=&-\unsurr\mx
\end{eqnarray*}
To finish, we compute
$\wh{C_{k,n}}.\,(r^2\,\wh{w_{24}}-r^3\,\wh{w_{14}}-r^3\,\wh{w_{23}}+r^4\,\wh{w_{13}})$,
when $5\leq k\leq n$ and we use the last relation of the third set
of relations above. The coefficient is given by
\begin{multline*}r^2\,ELOHH_{n-k,2k-6}\,-r^3\,ELOHH_{n-k,2k-5}\\-r^3\,ELOHH_{n-k,2k-5}\,+r^4\,ELOHH_{n-k,2k-4}\end{multline*}
and we see that it is indeed zero.\\
Let's deal with $(ii)$. First the fact that $\mz$ belongs to $K(4)$
can be achieved with Mathematica. Likewise, we check that $\my$
belongs to $K(5)$. Then, it remains to check that for all $j\geq 6$,
the algebra elements $\ovl{C_{1,j}}$, $\ovl{C_{2,j}}$,
$\ovl{C_{3,j}}$, $\ovl{C_{4,j}}$ and $\ovl{C_{5,j}}$ all annihilate
the vector $\my$. For $\wh{C_{1,j}}$ and $\wh{C_{2,j}}$, it follows
from IHNH. For $C_{3,j}$ and $\wh{C_{3,j}}$ we also use the first
and the last equations of the second set respectively. As for
$C_{5,j}$ and $\wh{C_{5,j}}$, we use the first and the third
relations of the first set respectively. Finally, we have using the
tables
\begin{eqnarray*}
[C_{4,j}.\,\my]_{w_{4,j}}&=&\unsur{r^2}\bigg(\frac{-1}{r^3.r^{j-6}}\bigg)-\frac{m}{r^{j-5}}\bigg(-\unsur{r^3}-\unsurr\bigg)+\unsur{r^{j-5}}\,\,=\,\,0\\
\left[
\wh{C_{4,j}}.\,\my\right]_{\wh{w_{4,j}}}&=&\unsur{r^2}\bigg(\unsur{r^{j-7}}\bigg)-\frac{m}{r^{j-5}}(r+r^3)-\frac{r^3}{r^{j-6}}\;\;=\;\;0
\end{eqnarray*}
Let's prove $(iii)$. Look at the action of the $\ovl{C_{s,t}}$'s on
$\mj_n$. First, when $s\geq 4$, we cut the sum term in $\mj_n$ into
three parts: a sum from $4$ to $s-1$, a term with the vectors ending
in node $s$ and a sum from $s+1$ to $n$. For the first part, we see
with ELOHH that the action is zero. For the middle part, we see with
the LO set that the action is zero. The last part however requires
more computations with the tangles. These are left to the reader. It
remains to check that the actions by $\ovl{C_{1,k}}$,
$\ovl{C_{2,k}}$, $\ovl{C_{3,k}}$ on $\mj_n$ are zero. Notice that
\begin{eqnarray*}
\ovl{C_{3,k}}&=&g_3\,\ovl{C_{2,k}}\,g_3^{-1}\\
g_3\,\mj_n&=&-\unsurr\,\mj_n\;\;\text{when $l=-\unsur{r^{2n-5}}$}
\end{eqnarray*}
Hence it suffices to study the actions by $\ovl{C_{1,k}}$ and
$\ovl{C_{2,k}}$. This is left to the reader.\\\\
This finishes the proof of Theorem $2$ point $(i)$. We now give an
example of how to compute the tangles of the table above. We show
below how ELOHH is computed. We want to compute
$\wh{C_{ij}}\wh{w_{i-s,i-t}}$.\vspace{1cm}
\begin{center}
\epsfig{file=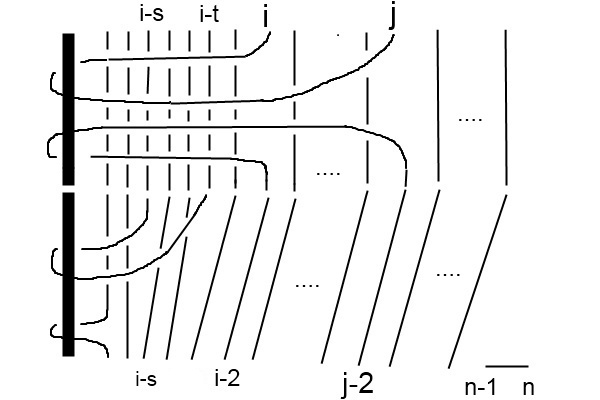, height=7cm}
\end{center}\vspace{1cm}
Use Reidemeister's move of type $3$ to move the bottom horizontal
strand of the upper tangle close to the pole as on the following
picture. Multiply the bottom of the tangle sucessively by the
products $g_{i-2}\dots g_3$ and $g_{j-2}\dots g_4$ at the cost of
divisions by $r^{i-4}$ and $r^{j-5}$ respectively. We must evaluate
\begin{center}
\epsfig{file=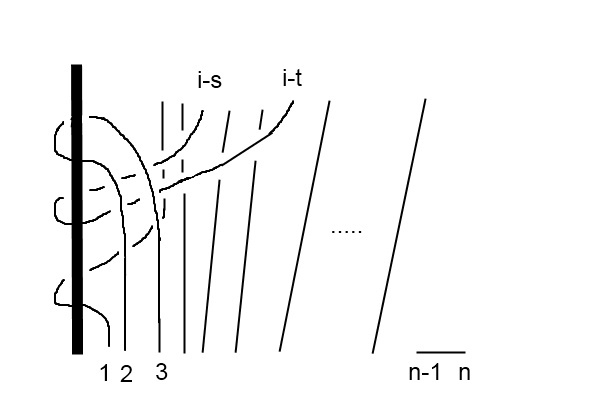, height=7cm}
\end{center}
where we omitted the top horizontal strand. Apply the commuting
relation in the upper left region of the figure and get
\begin{center}
\epsfig{file=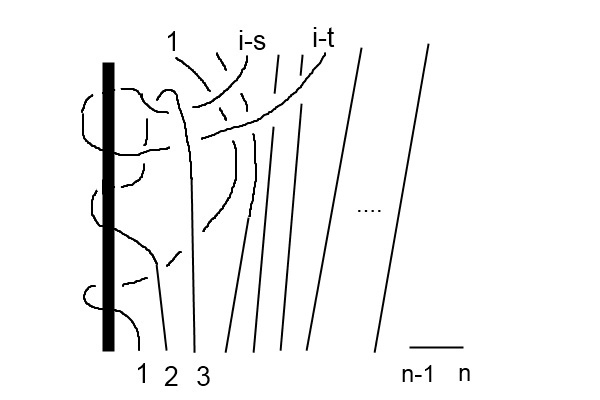, height=7cm}
\end{center}

\noindent Use the Kauffman skein relation twice, then use a
Reidemeister move of type $2$ twice, multiply at the bottom by $g_1$
at the cost of a division by $r$ and get a tangle that is zero. It
remains to compute the four terms arising from the two uses of the
Kauffman skein relation. When transforming the under-crossing into
an over-crossing on the upper left hand corner of the picture, one
must add two terms. The first term contains a loop around the pole
that can be "delooped" at the cost of a factor $\unsurr$. It is then
possible to apply the double twist relation. We hence obtain a
vertical strand joining nodes $i-t$ at the top and $2$ at the
bottom. We also obtain a loop that can be suppressed at the cost of
a factor $\unsur{l}$. The resulting strand is vertical and joins
nodes $i-s$ at the top and $3$ at the bottom. Use the sequence of
Reidemeister moves $R3$, $R2$, $R2$. Multiply at the bottom by $g_3$
at the cost of a division by $r$. This clarifies the first term.\\
When dealing with the second term, use a Reidemeister move of type
$3$ and get a loop around the pole. Suppress it at the cost of a
factor $r$. Further multiply at the bottom by $g_1$ at the cost of a
factor $\unsurr$ and apply the double twist relation twice, with a
Reidemeister move of type $2$ in between the two moves. Multiply at
the bottom by $g_3g_2$ at the cost of a division by $r^2$ and use a
Reidemeister move of type $2$ twice. Use a Reidemeister move of type
$3$. Multiply at the bottom by $g_3^{-1}$ at the cost of a
multiplication by $r$ and use a Reidemeister move of type $2$.
Multiply at the bottom by $g_1^{-1}$ at the cost of a multiplication
by $r$, use the double twist relation and a Reidemeister move of
type $2$. Multiply at the bottom by $g_3^{-1}$ at the cost of a
multiplication by $r$. Up to the coefficient, get the same tangle as
the one obtained after processing the first term. Gathering all the
moves that we did, we must now compute the expression of Figure $6$,
where we omitted the parts of the pictures that are not of direct
interest.
\begin{center}
\epsfig{file=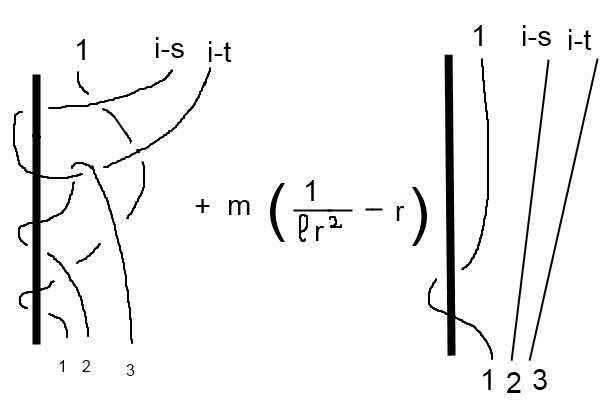, height=7cm}
\end{center}
\begin{center}
\textit{Figure $6$}
\end{center}
\noindent There is very little work that remains to be done on the
second tangle of Figure $6$. We must still straighten the vertical
strands. To that aim, multiply at the bottom sucessively by the
products $g_4^{-1}g_5^{-1}\dots g_{i-t}^{-1}$ and $g_3^{-1}\dots
g_{i-s}^{-1}$ at the cost of multiplications by $r^{i-t-3}$ and
$r^{i-s-2}$ respectively. \\
Finally, there are two more terms to compute. These arise from the
first tangle of Figure $6$ when we apply the Kauffman skein relation
for the second time. In the first term, there is a factor
$\unsur{l}$ arising from a loop. For the rest, after applying the
adequate moves, we get the same tangle as the one to the right. As
for the second term, there is a bit more work to be done. In what
follows, we use the abbreviation DT for double twist. The first step
is to do a Reidemeister move of type $3$. This then allows us to
apply the commuting relation. Multiply at the bottom by $g_2^{-1}$
at the cost of a multiplication by $r$, use $R2$, mutiply at the
bottom by $g_1$ at the cost of a division by $r$, do the sequence of
moves DT, R2, DT, multiply at the bottom by $g_3g_2$ and in order to
do so, divide by $r^2$, multiply at the bottom by $g_1^{-1}$ and in
order to do so multiply by $r$, use DT, then R2. After doing all
these moves, we get the tangle to the right of the picture. The
total is
$$m\,r^{2i-t-s-5}\bigg(\unsur{lr^2}-r+\unsur{l}-\unsurr\bigg)$$
One should not forget the factor $\unsur{r^{i+j-9}}$ from the
beginning. All together, it yields the coefficient of Proposition
$2$.

We end this section by showing that as a representation of the Artin
group, the representation $\n^{(n)}$ is equivalent to the
Cohen-Wales representation. Then, with the change of parameters
announced at the end of Theorem $1$, the point $(i)$ of Theorem $2$
implies the Main Theorem. In \cite{CGW}, the authors built all the
inequivalent irreducible representations of the quotient of ideals
$\mathcal{I}=C_ne_1C_n\,/<C_ne_ie_jC_n>_{i\not\sim j}$, where $C_n$
denotes the CGW algebra of type $D_n$. Only two of them have degree
the number of positive roots of a root system of type $D_n$, which
is also the degree of $\n^{(n)}$. By construction and by Theorem
$2$, point $(i)$, the representation $\n^{(n)}$ is an irreducible
representation of $\mathcal{I}$. Then, it must be equivalent to the
representation of \cite{CGW}. Our $r$ is the $\unsurr$ of
\cite{CGW}. Further, as a representation of the Artin group, the
represenation of \cite{CGW} is itself equivalent to the
representation of \cite{CW}, the one that was used to show the
linearity of the Artin group. The change of parameters is given in
the introduction of \cite{CGW} right before Theorem $1.1$.
\section{End of the proofs of the Theorems}
In this last section, we complete the proofs of Theorem $6$ and
Theorem $7$. These theorems provide of course important informations
about the structure of the Cohen-Wales representation of type $D_n$
and are extensively used in \cite{Clas}. In Theorem $6$, we must
still show that $S^{(0),(4,3)}$ and its conjugate, both of
dimensions $14$ cannot occur inside $V_7$. In Theorem $7$, we must
still show that the submodule of $V_n$ spanned by the
$\frac{n(n-1)}{2}$ vectors $t_{ij}$'s is irreducible. The latter
point uses the first point. We show the following results.
\begin{Proposition}
The Specht modules $S^{(0),(4,3)}$ and its conjugate
$S^{(0),(2^3,1)}$ don't occur in the Cohen-Wales space $V_7$.
\end{Proposition}
\begin{Proposition}
Proposition $3$ and Lemma $3$ imply Theorem $6$.
\end{Proposition}
\noindent \textsc{Proof.} By $\S\,3.2$, when $n\geq 6$, the only
irreducible $\mcalh(D_n)$-modules of degree $\cil$ are
$S^{(0),(n-2,2)}$ and its conjugate $S^{(0),(2,2,1^{n-4})}$, except
when $n=7$, when there are two more irreducibles, namely
$S^{(0),(4,3)}$ and $S^{(0),(2^3,1)}$. This settles Proposition $4$.
Let's prove Proposition $3$. Suppose there exists in $V_7$ an
irreducible invariant subspace $\W$ that is isomorphic to
$S^{(0),(4,3)}$. Then, we have $$\W\da_{\mcalh(D_5)}\simeq
2\,S^{(0),(3,2)}\oplus S^{(0),(4,1)}\qquad\;\;\;\;\;(\lozenge)$$ The
proof of Lemma $3$, point $(i)$ shows that if there exists an
invariant subspace of $\W\da_{\mcalh(D_5)}$ that is isomorphic to
$S^{(0),(3,2)}$, then it is unique. Hence it is
impossible to have $(\lozenge)$.\\
If now $\W$ is isomorphic to $S^{(0),(2^3,1)}$, then
$$\W\da_{\mcalh(D_5)}\simeq 2\,S^{(0),(2,2,1)}\oplus
S^{(0),(2,1,1,1)}\qquad\; (\ast)$$ But, by the proof of Lemma $3$,
point $(ii)$, the Specht module $S^{(2,2,1)}$ cannot be a
constituent of $\W\da_{\mcalh(D_5)}$. Thus, $(\ast)$ cannot happen
and
Proposition $3$ holds. This closes the proof of Theorem $6$.\\

Let's now show that the $\frac{n(n-1)}{2}$-dimensional invariant
subspace of Theorem $7$, say $\mathcal{T}$, is irreducible. Then it
must be unique. Excluding $n=4$, when $l=\unsurr$, the restrictions
on $r$ prevent the existence of an irreducible $d$-dimensional
invariant subspace of $V_n$ with $d\in\lb 1,n-1,n,\cil\rb$ by
Theorems $3,4,5,6$ respectively. So, if $\mathcal{T}$ has an
irreducible proper invariant subspace, the dimension of this
irreducible proper invariant subspace must be greater than or equal
to $\frac{(n-1)(n-2)}{2}=\cil+1$. But then it has a summand in
$\mathcal{T}$ whose dimension is less than or equal to $(n-1)$,
impossible. And so Theorem $7$ holds when $n\geq5$. \\

We conclude this paper by proving point $(ii)$ of Theorem $2$.
For the values of Theorem $2$, point $(i)$, the representation
$\n^{(n)}$ of $\cgwn$ that we built is reducible. Moreover, if
$\n^{(n)}$ were completely reducible, then by Proposition $1$, the
action of each $e_i$ on the Cohen-Wales space $V_n$ would be
trivial. This is impossible. Thus, $\cgwn$ is not semisimple for
these values of $l$ and $r$. As $r$ and $-\unsurr$ play identical
role, $\cgwn$ is not semisimple either for the values of Theorem $2$
point $(i)$ where $r$ has been replaced by $-\unsurr$.

\end{document}